\documentclass{article}
\usepackage{amsmath,amsthm,amsfonts,amscd,amssymb,eucal,latexsym}

\begin{document}

\baselineskip15pt
\textwidth=12truecm
\textheight=20truecm
\hoffset-.1cm
\voffset-.5cm

\def\Bbb{\mathbb}
\def\pp{\psi\vphantom A}
\def\ff{\varphi\vphantom A}
\def\ov{\overline}
\def\O{{\cal O}}
\def\var{\varepsilon}
\def\h#1{\hbox{\rm #1}}
\font\got=eufm10 scaled\magstephalf
\def\Ag{\h{\got A}}
\def\Fg{\h{\got F}}
\def\lam{\lambda}
\def\T{{\cal T}}
\def\C{{\cal C}}
\def\A{{\cal A}}
\def\R{{\cal R}}
\def\B{{\cal B}}
\def\M{{\cal M}}
\def\D{{\cal D}}
\def\N{{\cal N}}

\def\R{{\cal R}}
\def\P{{\Bbb P}}

\def\N{{\cal N}}
\def\C{{\cal C}}
\def\d{\hbox{\rm $\,$d}}
\def\End{\hbox{\rm End}}
\def\Aut{\hbox{\rm Aut}}
\def\O{{\cal O}}
\def\H{{\cal H}}
\def\I{{\cal I}}
\def\cC{{\cal C}}
\def\E{{\cal E}}
\def\cT{{\cal T}}
\def\N{{\cal N}}
\def\M{{\cal M}}
\def\D{{\cal D}}
\def\cB{{\cal B}}
\def\K{{\cal K}}
\def\A{{\cal A}}
\def\U{{\cal U}}
\def\Z{{\Bbb Z}}
\def\ov{\overline}
\def\Om{\Omega}
\def\var{\varepsilon}
\def\h#1{\hbox{\rm #1}}
\font\got=eufm10 scaled\magstephalf
\def\Ag{\h{\got A}}
\def\Fg{\h{\got F}}
\def\lam{\lambda}

\title{Regular Objects, Multiplicative Unitaries and Conjugation}
\author{
C. Pinzari \\
Dipartimento di Matematica, Universit\`a di Roma La Sapienza\\
00185--Roma, Italy\\ \\
 J.E. Roberts \\
Dipartimento di Matematica, Universit\`a di Roma Tor Vergata\\
00133--Roma, Italy}
\date{}
\maketitle

\begin{abstract}
 The notion of left (resp. right) regular object of a tensor
$C^*$--category equipped
with a faithful tensor functor into the category of Hilbert spaces
is introduced. If such a category has 
a left (resp. right) regular object, it can be interpreted as a
category
of corepresentations (resp. representations)
of some multiplicative unitary. A regular object is an object of the
category which is at the same time left and
right regular in a coherent way. A category with a regular object
is endowed with an associated 
standard braided symmetry.

Conjugation is discussed in the context of multiplicative unitaries and 
their associated Hopf
$C^*$--algebras. It is shown that 
the conjugate of a left regular object is a right regular object in the
same category. Furthermore the representation category of a locally
compact
quantum group has a conjugation. The associated
multiplicative unitary is a regular object in that category.

 \end{abstract}

\vfill

\noindent Research supported by MURST and CNR--GNAFA.\eject

\section{Introduction}
{}
  In this paper we look at the theory of multiplicative unitaries from 
the standpoint of their categories of representations and corepresentations. 
As is well known, multiplicative unitaries just 
express the fundamental property of the regular representation. Our 
approach therefore starts with a tensor category which may be thought
of as the tensor category of (unitary) representations of some quantum 
group. It is regarded as a concrete category in the sense that it is equipped 
with  a faithful tensor functor into the tensor category of Hilbert spaces. 
Once this tensor category has a ``regular object'' we will 
see that it allows an interpretation  as a category of 
representations of a multiplicative unitary and at the same time as a 
category of corepresentations of another multiplicative unitary. It is
instructive to compare this result with the Tannaka--Krein duality 
theorem or perhaps better with Woronowicz's duality theorem \cite{WTK}. 
In fact, our result starts with an embedded tensor category and constructs 
a multiplicative unitary and hence, if the multiplicative unitary is 
regular, two Hopf $C^*$--algebras \cite{BS}. However, by requiring the 
existence of a regular object, we are imposing a restriction 
that may not be easy to verify in practice and presupposes what a 
good duality theorem should prove. In fact, our result is close in spirit 
to Tatsuuma's duality theorem for locally compact groups \cite{T} where the group 
elements are identified in the regular representation using the 
multiplicative unitary.\smallskip
 
  Another aspect of the representation theory of multiplicative unitaries 
that has not received the attention it deserves is the conjugation
structure. 
We work here with multiplicative unitaries arising as the left regular 
representations of locally compact quantum groups. 
These are left regular objects in the category of corepresentations and we 
show that there is a canonical choice of conjugate which is a right
regular 
object. 
In fact, the multiplicative unitary
of a locally compact quantum group is a regular object in its
representation category.
Furthermore, we define the conjugate of any corepresentation up 
to unitary equivalence and the corresponding antilinear involution on 
intertwiners. 
This forms the subject matter of Section 5.
\smallskip

In this paper we prefer to work with strictly associative tensor products and a 
simple way of achieving this is to use as the underlying Hilbert spaces the 
Hilbert spaces in some fixed von Neumann algebra since these are objects 
in a strict tensor $W^*$--category. We will be concerned here with the 
representation categories of multiplicative unitaries and recall the 
basic definitions from \cite{BS}. If $K$ is such a Hilbert space then a
unitary 
$V$ on the tensor square $K^2$ is said to be multiplicative if 
$$V_{12}V_{13}V_{23}=V_{23}V_{12},$$
where we use the usual convention regarding indices and tensor products.
A representation of $V$ on a Hilbert space $H$ is a unitary $W\in(HK,HK)$ 
such that 
$$W_{12}W_{13}V_{23}=V_{23}W_{12},\quad \text{on}\quad HK^2.$$
If $W$ and $W'$ are representations of $V$ on $H$ and $H'$ respectively, we say 
that $T\in (H,H')$ intertwines $W$ and $W'$ and write $T\in(W,W')$ if 
$T\times 1_KW=W'T\times 1_K$. We define the tensor product of $W$ and $W'$ 
to be the representation $W\times W'$ on $HH'$ given by
$W\times W':=W_{13}W'_{23}$. The usual tensor product of intertwiners is 
again an intertwiner and in this way we get a strict tensor $W^*$--category 
${\cal R}(V)$ of representations of $V$. In fact this assertion does not 
depend on $V$ being multiplicative. When it is then $V$ itself is a 
representation of $V$ called the regular representation.\smallskip 

   A corepresentation of $V$ on $H$ is a unitary $W\in(KH,KH)$ such that 
$$V_{12}W_{13}W_{23}=W_{23}V_{12}\quad \text{on}\quad K^2H.$$ 
If $W$ and $W'$ are corepresentations on $H$ and $H'$ respectively, we say 
that $T\in(H,H')$ intertwines $W$ and $W'$ and write $T\in(W,W')$ if 
$1_K\times TW=W'1_K\times T$. The tensor product $W\times W'$ of 
corepresentations is defined by $W\times W':=W_{12}W'_{13}$. Just 
as in the case of representations we get a strict tensor $W^*$--category 
now denoted by $\C(V)$. If $\vartheta=\vartheta_{K,K}$ denotes the 
flip on $K^2$ then $\vartheta V^*\vartheta$ is again a multiplicative unitary 
and the mapping $W\mapsto \tilde W:=\vartheta_{H,K}W^*\vartheta_{K,H}$ defines  
a 1--1 correspondence between representations of $V$ and corepresentations 
of $\vartheta V^*\vartheta$. However, it does not define an  isomorphism 
of tensor $W^*$--categories since $W\times W'\mapsto \tilde W'_{13}\tilde W_{12}$ 
and so leads to an alternative definition of the tensor product of 
corepresentations. In fact the two expressions for the tensor product 
will be equal if and only if $\vartheta_{W,W'}\in(W\times W',W'\times W)$, 
cf. Prop. 2.5 in \cite{W}.\smallskip 

\section{Regular Objects and Multiplicative Unitaries}

The main aim of this section will be to provide characterizations of categories 
of representations and corepresentations of multiplicative unitaries and in 
particular to study tensor categories which are simultaneously a tensor category
of representations of a multiplicative unitary and of corepresentations of 
some (other) multiplicative unitary. The main idea is to replace the notion 
of multiplicative unitary by that of regular object. Thus multiplicative 
unitaries are seen as intertwining operators, taking us back to the origins 
of the theory. They are therefore seen not as determining a category of 
representations or corepresentations but as being a structural element in some 
tensor category. This helps us to understand the degree to which they are 
not unique and to see the tensor categories that are simultaneously a 
category of representations and a category of corepresentations as being 
tensor categories with a left and right regular object.\smallskip 
 

   Here is our motivating example. Let ${\cal H}$ denote
the strict tensor $W^*$--category of Hilbert spaces  in a von Neumann algebra $M$ 
and $\vartheta$ its unique permutation symmetry. Let $K$ be an object of 
${\cal H}$ and $V$ a multiplicative unitary on $K^2$. 
We let  ${\cal R}(V)$  and ${\cal C}(V)$  be the tensor
$W^*$--categories of representations and corepresentations of $V$  on
Hilbert spaces of
 ${\cal H}$. These are to be considered as equipped with the forgetful functor $\iota$ into 
 ${\cal H}$ itself, regarded as the subcategory of trivial representations or 
corepresentations. Thus $\iota$ is an idempotent tensor $^*$--functor.\smallskip

   We now ask the following question: when can a strict tensor $W^*$--category 
${\cal T}$ equipped with a faithful idempotent tensor $^*$--functor 
$\iota_{\cal T}:=\iota$ onto a tensor $W^*$--subcategory of Hilbert spaces be 
interpreted as a category of representations or corepresentations of a 
multiplicative unitary? Note that $\vartheta_{H,H'}$ is an intertwining 
operator in ${\cal R}(V)$  
for a tensor product of representations $W$ on $H$ and $W'$ on $H'$ whenever 
either $W$ or $W'$ is a trivial representation. 
Any full tensor subcategory $\cal{T}$ of ${\cal R}(V)$ containing the regular 
representation $V$ has the striking property that for any object $W$, 
$W\times V$ is a, possibly infinite, direct sum of copies of $V$ but more is 
true: we set $\eta_W:=W\in(W\times V,\iota(W)\times V)$ then $\eta\in(R,R\iota)$, 
where $R$ denotes the functor of tensoring on the right by $V$, is a natural 
unitary transformation such that 
$$\eta_{W\times W'}=(\eta_W)_{13}\circ 1_W\times\eta_{W'},\eqno(2.1)$$ 
for each pair $W,W'$ of objects of ${\cal T}$.\smallskip 

      To formalize the essential aspects of the above situation we consider a 
strict tensor $W^*$--category $\cal{T}$ equipped with a faithful
idempotent tensor
$^*$--functor $\iota_{\cal{T}}:=\iota$. The tensor subcategory
$\iota(\cal{T})$ is 
equipped with a (permutation) symmetry $\vartheta$. We further suppose 
that given objects $W$ and $W'$ of $\cal{T}$ there are arrows 
$\vartheta_{W,\iota(W')}\in(W\times\iota(W'),\iota(W')\times W)$ 
and $\vartheta_{\iota(W),W'}\in(\iota(W)\times W',W'\times\iota(W))$, 
necessarily unique,  whose image under $\iota$ is
$\vartheta_{\iota(W),\iota(W')}$. 
We call a {\it right regular} object of $\cal{T}$ a pair $(V,\eta)$ 
consisting of an object $V$ of $\cal{T}$ and 
a unitary natural transformation $\eta\in(R,R\iota)$, where $R$ denotes 
the functor of tensoring on the right by $V$, satisfying (2.1) above    
for each pair $W,W'$ of objects of ${\cal T}$. Here $(\eta_W)_{13}$ is to
be 
understood as $1_{\iota(W)}\times\vartheta_{V,\iota(W')}\circ\eta_W\times
1_ {\iota(W')}
\circ 1_W\times\vartheta_{\iota(W'),V}$.
$(2.1)$ implies that $\eta$ evaluated on the tensor unit ${\mathbb C}$
is 
$1_{\mathbb C}$.\smallskip 

   The following result now provides an answer to the above 
question.\smallskip

\noindent
{\bf 2.1 Theorem} {\sl Any tensor $W^*$--category equipped with a 
faithful idempotent tensor $^*$--functor into a tensor subcategory of 
Hilbert spaces and with a right regular object is isomorphic to 
a tensor $^*$--subcategory of ${\cal R}(V)$ for some multiplicative unitary $V$.}\smallskip 

\noindent{\bf Proof.} Let $\eta$ denote a natural transformation in $(R,R\iota)$ 
making an object $V$ into a right regular object then $\iota(\eta_W)$ 
is a unitary for each object $W$. If $T\in(W,W')$ then the naturality of $\eta$ 
shows that $\iota(T)$ intertwines $\iota(\eta_W)$ and $\iota(\eta_{W'})$, 
once we know that these are representations of $\iota(\eta_V)$. In 
particular, taking $\eta_W$ as $T$, naturality gives
$$\iota(\eta_W)\times 1_V\circ\eta_{W\times V}=
\eta_{\iota(W)\times V}\circ\eta_W\times 1_V.\eqno(2.2)$$ 
Equation $(2.1)$ tells us that the tensor product in the category
corresponds to 
the tensor product of representations. Bearing this in mind, $(2.2)$ tells
us 
that $\iota(\eta_W)$ is a representation of $\iota(\eta_V)$ and the particular case $W=V$ 
tells us that $\iota(\eta_V)$ is indeed a multiplicative unitary.
\medskip

   The notion of multiplicative unitary and Theorem 2.1 can be easily 
generalized replacing a multiplicative unitary on a Hilbert space by 
a multiplicative invertible in a monoidal category. We refrain from
spelling 
this out to keep a uniform setting for this paper.\smallskip

   Notice that the isomorphism in question is even canonical, given
$\eta$, 
and commutes with $\iota$. 
But there are several other comments to be made about this result. First,
$(2.2)$ has the 
structure of an associative law: it equates two ways of passing from $RR$
to
$R\iota R\iota$. Secondly,  there is an analogous result 
for corepresentations. We define a notion of left regular object 
by dualizing in ${\cal T}$ with respect to the composition law $\times$.
Our 
category ${\cal T}$ is isomorphic to a tensor subcategory of the category 
of corepresentations of a multiplicative unitary if it admits a  
left regular object. If $\xi$ denotes a natural transformation rendering 
$V$ a left regular object $(\xi,V)$, then the unitary associated with an 
object $W$ is $\iota(\xi_W^{-1})$. The appearance of an inverse here is
just an 
artefact of conventions.\smallskip 

    One might have thought of basing a definition of right regular object 
on a different familiar property of the regular representation, namely 
that $W\times V$ is a, possibly infinite, direct sum of copies of $V$ 
for each object $W$ of $\cal{T}$. This property is too weak in that 
it does not imply the coherence properties of the previous definition 
and furthermore puts an unwanted emphasis on the notion of infinite 
direct sum. In fact, our definition implies the second property once 
we specify, as we now do, that $\iota(\cal{T})$ is just a category of
Hilbert 
spaces, i.e.\ a (strict) tensor $W^*$--category with unit reducing to 
the complex numbers where every object is a (possibly infinite) 
direct sum of the unit. Since $\iota(W)$ is a direct sum of copies of the 
unit ${\mathbb C}$, 
$W\times V\simeq\iota(W)\times V$ is a direct sum of copies of $V$. 
Note that if $V_r$ is right regular and $V_\ell$ is left regular then 
$V_\ell\times V_r$ is a direct sum of copies of both $V_\ell$ and $V_r$.
It 
follows that, if we have both a left and a right regular object, then
these 
objects are unique up to quasiequivalence in the $W^*$--category in
question. 
In a $\sigma$--finite $W^*$--category a left or right regular object with 
infinite multiplicity will then be unique up to unitary
equivalence.\smallskip

The following variant on the definition of a right regular object is
worth 
noting. Consider a tensor $W^*$--category ${\cal T}$ and unit ${\mathbb
C}$, 
but where the endofunctor $\iota$ is not a priori defined. Suppose for 
each object $W$, there is an unitary arrow 
$\eta_W\in(W\times V,\iota(W)\times V)$ such that $\iota(W)$ is a
(possibly
infinite) direct sum of the tensor unit. Suppose (2.1) holds and 
$$\eta_{W'}\circ T\times 1_V\circ\eta_W^{-1}\in(\iota(W'),\iota(W))\times
1_V,
\quad T\in(W,W').$$
Then setting $\iota(T):=\eta_{W'}\circ T\times 1_V\circ\eta_W^{-1}$, 
we get a tensor $^*$--endofunctor from ${\cal T}$ into a tensor
subcategory 
of Hilbert spaces. If $\eta_{\iota(W)}=1_{\iota(W)\times V}$ for each
object 
$W$ of ${\cal T}$, $\iota$ is even idempotent. This illustrates the role
of 
(2.1) in guaranteeing that a tensor $W^*$-category can be embedded into a
tensor 
category of Hilbert spaces.\smallskip

   We now make some further remarks on the notion of regular object, supposing 
for the moment that our category ${\cal T}$ has sufficient irreducibles in the 
sense that every object is a (possibly infinite) direct sum of irreducibles and
${\cal T}$ 
is closed under finite direct sums. Suppose further that the full subcategory 
${\cal T}_f$ whose objects are finite direct sums of irreducibles is a tensor 
subcategory and that $\iota_W$ is finite dimensional for each irreducible $W$.  
A {\it dimension function} $d$ on ${\cal T}_f$ assigns to each object $W$ of 
${\cal T}_f$ a  $d(W)\in{\mathbb R}_+$ such that 
$$d(W\oplus W')=d(W)+d(W'),$$
$$d(W\times W')=d(W)d(W').$$ 
Note that if $F:{\cal T}_f\to{\cal T'}_f$ is a tensor $^*$--functor and $d'$ is a 
dimension function on ${\cal T'}_f$, then $F\circ d'$ is a dimension
function on ${\cal T}_f$.
Thus, our category ${\cal T}_f$ has an integer--valued dimension function induced by the 
tensor $^*$--functor $\iota$ from the Hilbert space dimensions. 
If ${\cal T}_f$ has conjugates, then there is another dimension function, 
not necessarily integer--valued, given intrinsically by its structure as a 
tensor $C^*$--category \cite{LR}. 
Let $I$ be an index set labelling the equivalence classes of irreducibles and $W_i$ 
an irreducible of class $i\in I$.
Then a simple computation shows that if $d$ is a dimension function and 
$d_i:=d(W_i)$ then 
$$d_jd_k=\sum_i\,m^j_{ki}d_i,\quad i,j,k\in I,$$  
where $m^j_{ki}$ denotes the dimension of $(W_i,W_j\times W_k)$. Thus 
the dimension function, which is determined by the $d_i$, $i\in I$, gives an 
eigenvector with positive entries of the matrix $m^j$ corresponding to the 
eigenvalue $d_j$ and simultaneously an eigenvector of $m_k$ with eigenvalue 
$d_k$. Conversely, any such simultaneous eigenvalue does arise in this way. 
Suppose that $V$ is a left regular object of ${\cal T}$ such that $(W,V)$ is finite 
dimensional for each irreducible $W$ and hence for each object $W$ of ${\cal T}_f$. 
 Let $v_i$ be the dimension of $(W_i,V)$ and let 
$d(W)$ be defined so that $V\times W$ is a direct sum of $d(W)$ copies of $V$. 
Then $d$ is an integer-valued dimension function and 
$$d_jv_i=\sum_k m^k_{ji}v_k.$$ 
The case of a right regular object can be treated similarly. 
If ${\cal T}_f$ has conjugates then we have a corresponding involution 
$i\mapsto\ov i$ on $I$ and 
$$m^i_{jk}=m^{\ov i}_{kj}=m^k_{\ov ji}.$$ 
If $d$ is a dimension function, there is a conjugate dimension function
$\ov d$ 
such that $\ov d(W)=d(\ov W)$ for each object $W$ of ${\cal T}_f$. The 
interesting dimension functions, such as the intrinsic dimension 
function of a tensor $C^*$--category with conjugates\cite{LR}, are
self-conjugate.\smallskip
 
   Now, there is another natural transformation implicitly involved 
in $(2.1)$, namely, $\theta\in (R\iota,L\iota)$, defined by 
$$\theta_W:= \vartheta_{\iota(W),V}.$$
This brings us to the concept of braided symmetry, developed in the Appendix 
of \cite{DPR}. 
Let $\var$ be a braided symmetry relative to a left regular object $V$ of 
${\cal T}$. Thus $\var$ is a unitary natural transformation
from the functor $R$ of tensoring on the right by $V$ to the functor 
$L$ of tensoring on the left by $V$
such that $$\var_{W\times W'}=\var_W\times 1_{W'}\circ
1_W\times\var_{W'}.$$ 
Note that
$\var_{\iota(W)}=\vartheta_{\iota(W),V}=\theta_W$. 
Since $V\times W$ is just a multiple of $V$ and the functor $L$ is 
faithful, $\var_W$ is uniquely determined by $\var_V$ using
$$\var_{V\times W}=\var_V\times 1_W\circ 1_V\times\var_W.$$
The index notation for tensor products will now be taken to refer to the 
braided symmetry. This is consistent with its use in $(2.1)$. 
Using the braided symmetry, we get
 a unitary natural  transformation 
$\eta$ from $R$ to $R\iota$ defined by
$$\eta_W=\var^{-1}_{\iota(W)}\xi_W\var_W,$$ 
where $\xi$ is the unitary natural transformation from $L$ to $L\iota$ 
making $V$ into a left regular object.  
We ask whether $\eta$ makes $V$ into a right regular object. 
This question is addressed in the following results.\medskip

\noindent{\bf 2.2 Proposition} {\sl Let $(\xi,V)$ be a left regular object
of
${\cal T}$. The braided symmetries $\var$ for 
${\cal T}$ relative to $V$  are in $1-1$ correspondence with invertible 
natural transformations $\eta$ from $R$ to $R\iota$ such that
$$\eta_{W\times W'}=
(\xi_{W'})_{32}(\eta_W)_{13}(\xi^{-1}_{W'})_{32}(\eta_{W'})_{23}\,.\eqno(2.3)$$
$\var$ and $\eta$ are related by 
$$\var_W=\xi^{-1}_W\theta_W\eta_W\,.\eqno(2.4)$$}

\noindent{\bf Proof}. Given $\var$, equation $(2.4)$ defines a natural
unitary 
transformation $\eta$ and $(2.3)$ follows by direct
computation. Conversely, given $\eta$, equation $(2.4)$ defines a  
natural unitary transformation $\var$ which is a braided symmetry by 
virtue of $(2.3)$.\medskip

If we take the images under $\iota$ of the terms in $(2.3)$ then the  
computation leading to $(2.3)$ can be modified to show that the analogous 
identity holds with the tensor product notation now referring to the 
permutation of Hilbert spaces. Note, too, that $(2.3)$ can be used to 
compute $\eta$ in terms of $\eta_V$.\medskip

\noindent{\bf 2.3 Theorem} {\sl Given functors $L$ and $R$ of tensoring on  
the left and right, respectively, by an 
object $V$ of ${\cal T}$ and invertible natural transformations $\xi \in
(L,L\iota)$,
$\eta\in (R,R\iota)$ and $\var\in (R,L)$ such that $\xi\var=\theta\eta$, 
consider the following four conditions:
\begin{description}
\item{a)} $(\xi,V)$ is a left regular object,
\item{b)} $(V,\eta)$ is a right regular object, 
\item{c)} $\var$ is a braided symmetry relative to $V$,
\item{d)} $\eta_W\times 1_{\iota(W')}\circ 1_W\times
\xi_{W'}=1_{\iota(W)}\times 
\xi_{W'}\circ\eta_W\times 1_{W'}$, for each pair $W,W'$ of objects of
$\cal T$.
\end{description}
Then any three of these conditions imply the fourth.}\smallskip

\noindent{\bf Proof.} We see from Proposition 2.2 that given a) and c), 
b) is equivalent to  
requiring that each pair $(\eta_W)_{13}$ and $(\xi_{W'})_{32}$ commute. 
Interchanging 2 and 3 using the braided symmetry, we see that, given 
a) and c), b) and d) are equivalent. Similarly, given b) and c), a) and d) 
are equivalent. It remains to show that a), b) and d) imply c). However, 
given a), b) and d), $(2.3)$ follows from $(2.1)$, since d) implies 
that $(\eta_W)_{13}$ and $(\xi_{W'})_{32}$ commute. Thus $\var$ is 
a braided symmetry, completing the proof. 
\medskip

   An alternative way of proving the above theorem is by arguing in 
terms of a  diagram with ten vertices, where the conditions a), b), 
c) and d) are expressed as the commutativity of subdiagrams. The 
reader is urged to draw the diagram for himself. Begin with an outer 
square whose sides are used as an hypotenuese for the conditions on 
$\eta$, $\xi$, $\var$ and $\theta$ respectively with d) as a rhombus in 
the middle of the square.  

   We may also strengthen one of the implications in the above theorem.\medskip

\noindent{\bf 2.4 Lemma} {\sl Under the hypotheses of Theorem 2.3, the 
conditions a), b) and 
\begin{description} 
\item{d$'$)} $\eta_V\times1_{\iota(V)}\circ 1_V\times\xi_V=1_{\iota(V)}
\times\xi_V\circ\eta_V\times 1_V$
\end{description}
imply that $\varepsilon$ is a braided symmetry.}\medskip

  The necessary computations can be found in the proof of 
Theorem A.2 in \cite{DPR}. This proof can be rewritten entirely 
in terms of compositions 
in the tensor category ${\cal T}$ and this is recommended to the 
reader as an exercise. The computations also show that a), c) 
and d) with $V$ in place of $W'$ imply b) and that b), c) and d) with 
$V$ in place of $W$ imply a).\smallskip 

  We will refer to a braided symmetry fulfilling the conditions of 
Theorem 2.3 as being a {\it standard} braided symmetry. In the presence 
of a standard braided symmetry we have an object which is at the 
same time a left and right regular object in a coherent way in that it 
fulfills d) of Theorem 2.3. We call such an object a {\it regular} 
object. Under these circumstances 
we have the following corollary of Theorem 2.1.\smallskip 

\noindent
{\bf 2.5 Corollary} {\sl A tensor $W^*$--category equipped with a faithful 
idempotent tensor $^*$--functor into a tensor subcategory of Hilbert
spaces and 
with a regular object $(\xi,V,\eta)$ is isomorphic to a 
tensor subcategory of both ${\cal R}(\iota(\eta_V))$ and 
${\cal C}(\iota(\xi_V)^{-1})$. 
There is an associated standard braided symmetry given by
$(2.4)$.}\smallskip 
   
  In particular, if a multiplicative unitary $V$ considered as a 
corepresentation and hence a left regular object in ${\cal C}(V)$ 
is even a right regular object, there is a multiplicative 
unitary $\hat V$ on the same Hilbert space such that ${\cal C}(V)$ 
is canonically isomorphic as a tensor $W^*$--category to a tensor 
$^*$--subcategory of ${\cal R}(\hat V)$. We do not know when the 
image coincides with ${\cal R}(\hat V)$.\smallskip 

  Given a left regular object $V$ of $\cal{T}$, we would like
to analyse in how many different ways we may 
choose $\eta$ and $\var$ so as to fulfill the conditions of Theorem 2.3. 
Of course $\var$ determines $\eta$ uniquely so our question amounts to
parametrizing the standard braided 
symmetries relative to  $V$. It is convenient to rephrase this problem in terms
of the associated natural transformations (Prop. 2.2). 
In the remark following Proposition 2.2, we have noted that any such
natural 
transformation is uniquely determined
by its value in $V$. Therefore we first consider the situation where 
$\cal{T}$   is ${\cal C}(V)_V$, the full tensor subcategory of 
${\cal C}(V)$ generated by $V$. 
\medskip

\noindent{\bf 2.6 Lemma} {\sl Let   $\xi$ be the 
 natural unitary transformation from $L$ to $L\iota$ 
making $V$ into a left regular object of ${\cal C}(V)_V$, and let 
$\eta_V$ be a 
unitary operator on the Hilbert space $K^2$ of 
$V$. Then there is a  natural  unitary transformation
$\eta$ from $R$ to $R\iota$ taking the value $\eta_V$ at $V$ 
 and defining a standard braided symmetry $\varepsilon$ on  
${{\cal C}(V)}_V$
if and only if $\eta_V$  satisfies
\begin{description}
\item{\rm a)} $\eta_V \times 1_{\iota(V)}\circ 1_V \times  \xi_V=
1_{\iota(V)} \times \xi_V\circ
\eta_V  \times 1_{V}\ ,$
\item{\rm b)} $\eta_V\circ T\times 1_V=\iota(T)\times
1_V\circ\eta_V,\quad T\in(V, V)$, 
\item{\rm c)}  $\eta_V\in(V^{\times 2}, \iota(V)\times V)$, 
\item{\rm d)} $\iota(\xi_V)\times 
1_V\circ(\eta_V)_{13}\circ1_V\times\eta_V=
(\eta_V)_{13}\circ\xi_V\times 1_V$.
\end{description}}

\noindent{\bf Proof.} If  $\eta$ is the natural  unitary transformation associated
with the standard braided symmetry $\varepsilon$ then, by Proposition 2.2, 
$\eta\in(R, R\iota)$, thus, evaluating in $V$, we get c); a) is a special case
of d) in Theorem 2.3. Since $\eta$ is natural, given any pair 
$W$, $W'\in{\cal C}(V)_V$ and
any  $T\in(W, W')$, $\iota(T)\times 1_V\circ\eta_W=\eta_{W'}\circ T\times1_V$,
thus choosing $W=W'=V$ we obtain b). On the other hand 
 $V\in(V\times\iota(V), V^{\times
2})$, therefore, as $\eta$ makes $V$ into a right 
regular object, d) follows from the naturality of $\eta$.  Conversely, 
by virtue of a) and Lemma 2.4, it suffices to show that
$\eta_{V^{\times r}}:=(\eta_{V})_{1 r+1}\dots(\eta_{V})_{r r+1}$
is a natural unitary transformation from $R$ to $R\iota$ making $V$
into a right regular object.
It is easy to see that c) yields $\eta_{V^{\times r}}\in(V^{\times r}\times V, 
\iota(V^{\times r})\times V)$, 
and we must show the naturality
of $\eta$,
 i.e. that for any $T\in(V^{\times r}, V^{\times s})$, 
$\eta_{V^{\times s}}\circ T\times 1_V=
\iota(T)\times 1_V\circ\eta_{V^{\times r}}$.  We note that if this 
relation  holds for $r+1$ and $s+1$ then it holds for $r$ and $s$ as
well, since  $(V^{\times r}, V^{\times s})$ embeds in $(V^{\times r+1}, V^{\times
s+1})$ via $R$. Therefore, it suffices to assume $r, s$ sufficiently large. Now
by b) the  relation holds for $r=s=1$.  
We regard the Hilbert space $K$ of the corepresentation $V$ as a space 
of bounded linear operators from $K$ to $K^2$ by letting the elements 
of $K$ act by tensoring on the left. By the multiplicativity
of $V$, $V\vartheta_{K, K} K$ is a Hilbert space $\tilde{K}$ of intertwiners
of ${\cal C}(V)_V$ contained in $(V, V^{\times 2})$ and property d) shows that
the  desired relation holds for elements of  $\tilde{K}$.
For $s\geq 2$, on the other hand, $(V^{\times r},  V^{\times s})$ is
generated as a weakly closed subspace of $(K^r, K^s)$ by
elements of the form 
$\psi\times 1_{V^{s-2}}\circ T$,
with $T\in(V^{\times r},  V^{\times s-1})$ and $\psi\in\tilde{K}$. The 
relations therefore hold for a  generating set
of intertwiners in ${\cal C}(V)_V$, and hence  for all the intertwiners, 
completing  the proof.\medskip

It emerges from  the proof that property d) has the role of ensuring 
the naturality of $\eta$  for elements of the Hilbert space 
$\tilde{K}\subset (V,V^{\times 2})$.
$\tilde{K}$ could  be replaced by any other
Hilbert space with support $I$ in $(V, V^{\times 2})$. Choosing  
$\hat{K}:=\eta_V^{-1}K$ amounts to replacing d) by

$d')$  $\eta_V$ is a multiplicative unitary on $K^2$.

\noindent
Furthermore, we note that  b) characterizes the elements of 
$(V, V)\subset (K, K)$. Indeed, if 
$T\in(K, K)$ satisfies $\eta_V\circ T\times 1_V=T\times 1_V\circ\eta_V$
then by c) $T\times 1_V={\eta_V}^*\circ T\times 1_V\circ{\eta_V}\in
(V^{\times 2}, V^{\times 2})$, so $T\in(V, V)$.

   To parametrize the  standard braided symmetries, we shall need two 
further notions: let $W$ be an object of ${\cal C}(V)$ acting  on $H$ and set
$$G_W:=\{U\in{\cal U}(H) : TU^{\times r}=U^{\times s}T\ ,\quad
T\in(W^{\times r}, W^{\times s})\}\ .$$

\noindent{\bf 2.7 Proposition} {\sl If $V$ is the regular corepresentation
then 
$G_V=\{U\in(V,V)' : \delta(U):=VU\times 1_VV^{-1}=U\times U\}\ .$}\smallskip

\noindent{\bf Proof.} $U\in G_V$ implies $\delta(U)=U\times U$ 
 since $V\in (V\times\iota_V,V\times V)$ and $U\in(V,V)'$. 
For the converse, note that the above two conditions suffice to conclude 
that $U\in G_V$ by the fundamental 
property of the regular corepresentation.\medskip
 
We regard the dual multiplicative unitary 
$V^d:=\vartheta_{K, K}V^{-1}\vartheta_{K, K}$ as a left regular object
of ${\cal C}(V^d)$. We call $V$ {\it weakly
irreducible\/} if $(V, V)\cap(V^d, V^d)={\mathbb {C}}I$. If $V$ is irreducible
in the sense of \cite{BS} then it is weakly irreducible.
\medskip

\noindent{\bf 2.8 Theorem} {\sl Let $V$  be a weakly irreducible
multiplicative unitary, and let
 $\eta$,  $\tilde{\eta}\in(R, R\iota)$
 be  natural unitary transformations defining  standard braided
symmetries  on ${\cal C}(V)_V$. Then there is a unique unitary
$U\in G_V$ such that 
${\tilde{\eta}}_V=U\times 1_V\circ\eta_V$. Conversely, given any unitary 
$U\in G_V$ and any natural transformation
$\eta$ as above there is a unique  
$\tilde{\eta}$ such
that ${\tilde{\eta}}_V=U\times 1_V\circ\eta_V$.}
\smallskip

\noindent{\bf Proof.} Let $\eta$ and $\tilde{\eta}$ define
standard braided symmetries on $\C(V)_V$.
By virtue of the commutation
relation a)  of  Lemma 2.6, for any $\omega\in(K, K)_*$,
$\omega\otimes\iota(\eta_V)\circ V=
V\circ\omega\otimes\iota(\eta_V)$, therefore  
$\omega\otimes\iota(\eta_V)\in(V^d,V^d)$. On the other hand, by
c) of the same lemma, 
 ${\tilde{\eta}}_V{\eta_V}^{-1}\in(\iota(V)\times V,
\iota(V)\times V)=(\iota_V, \iota_V)
\otimes(V, V)$, so ${\tilde{\eta}}_V{\eta_V}^{-1}\in
(K, K)\otimes(V, V)\cap(V^d, V^d)=(K, K)\otimes{\mathbb C}I$ 
since $V$ is weakly irreducible. Now, by  b), 
$\iota\otimes\omega(\eta_V)\in(V, V)'$, so both $\eta_V$, $\tilde{\eta}_V$,
and therefore
$\tilde{\eta}_V{\eta_V}^{-1}$ belong to $(V, V)'\otimes(K, K)$.
We conclude that there is a $U\in(V, V)'$ with $\tilde{\eta}_V=
U\times 1_V\eta_V$. Finally, comparing  d) for
$\eta_V$ and $\tilde{\eta}_V$, we conclude that $U$
satifies $V\circ U\times1_V=U\times U\circ V$,
i.e. $U\in G_V$.

Conversely, a straightforward computation shows that  any  unitary
of the form $\tilde{\eta}_V:=U\times 1_V\eta_{V}$, with $U\in G_V$,
satisfies the properties stated in the previous lemma.
\medskip 

\noindent{\bf 2.9 Lemma} {\sl Let $V$ be a multiplicative 
unitary and $\eta_V$ a unitary operator 
on $K^2$ satisfying properties a) and c) of Lemma 2.6. If $V$ is
regular in the sense of \cite{BS} then the algebra 
generated by 
$(V, V)$ and $(V^d,V^d)$ acts irreducibly on $K$.}\smallskip

\noindent{\bf Proof.} For any pair $\psi$ and $\varphi$ of elements 
of $K$   we may write 
$$\psi^*\times1_VV\vartheta\varphi\times1_V
=\sum_i\psi^*\times1_V{\eta_V}^*\varphi_i\times1_V
{\varphi_i}^*\times1_V\eta_V V\vartheta\varphi\times1_V\ ,$$
where $\varphi_i$ is an orthonormal basis
of $K\ .$ Now 
$$K^*\times1_V{\eta_V }^*K\times1_V\subseteq(V^d, V^d)$$
as ${\eta_V}_{12}$ and $V_{23}$ commute.
 On the other hand
$K^*\times1_V\eta_VV\vartheta K\times1_V\subseteq(V, V)$,
therefore $K^*\times1_VV\vartheta K\times1_V$ is a subspace of 
the weak closure of the algebra generated by $(V, V)$ and 
$(V^d,V^d)$. On the other hand 
this subspace generates the compact operators since $V$ is regular,
completing the proof.
\medskip

\noindent{\bf 2.10 Theorem} {\sl If $V$ is a  multiplicative unitary
and the algebra generated by $(V, V)$ and $(V^d,V^d)$
acts irreducibly on $K$ then
any  braided symmetry $\var$  on ${\cal C}(V)_V$ extends uniquely
to a  braided symmetry  on ${\cal C}(V)\ ,$ 
standard if $\var$ is standard.}\smallskip

\noindent{\bf Proof.} 
 We have already noted that any braided symmetry on 
$\C(V)$ is  determined uniquely
by $\eta_V$. The explict relation is
$$1_V\times\eta_W=
(\xi_W)_{32}({{\eta_V}^{-1}})_{13}({\xi_W}^{-1})_{32}({\xi_W}^{-1})_{12}({\eta_V})_{13}
({\xi_W})_{12}.$$ We show that the right hand side does define a natural 
unitary transformation from $R$ to $R\iota$ in $\C(V)$ satisfying
$(2.3)$. For brevity, we write $W^{-1}$ for $\xi_W\ .$

Let $H$ be the Hilbert space of the corepresentation $W\ .$
The key idea is to show
that the unitary operator $X_W$ on $KHK$ defined by the right hand side
acts trivially on the first factor, by showing that  
its first component lies in the commutant of the algebra generated
by $(V, V)$ and 
$(V^d, V^d)$,
this being the
complex numbers, by assumption. We first show that
the first component of $X_W$
is in the commutant of $(V^d, V^d)$.
 Since the first component of $W$ is in $(V^d,V^d)$,
 and ${W_{32}}^{-1}$ acts trivially on the first factor,
it is enough to prove the claim for the first component of 
$({\eta_V}^{-1})_{13}W_{32}W_{12}({\eta_V})_{13}\ .$ However, 
by the corepresentation
relation  $W_{32}W_{12}={V_{31}}^{-1}W_{12}V_{31}\ ,$ 
we are thus reduced to showing that the first component of 
$({{\eta_V}^{-1}})_{13}{V_{31}}^{-1}W_{12}V_{31}({\eta_V})_{13}$ is 
contained in the commutant of 
$(V^d, V^d)$.
Now, this holds in the special case $W=V$ since this operator coincides with
$V_{32}({\eta_V})_{23}V_{12}\ .$ 
In the general case, since the first component of $W$ is contained in
$(V^d, V^d)'$ we deduce that we can approximate $W$ weakly by finite
sums of operators of the form $1_K\times A^*V1_K\times B\ ,$ with $A, B\in(H, K)\ .$
Hence it suffices  if the first component of 
$$({\eta_V}^{-1})_{13}{V_{31}}^{-1}1_K\times A^*\times 1_KV_{12}1_K\times B\times 1_KV_{31}({\eta_V})_{13}=$$
$$1_K\times A^*\times 1_K({\eta_V}^{-1})_{13}{V_{31}}^{-1}V_{12}V_{31}({\eta_V})_{13}1_K\times B\times 1_K\ ,$$
is in 
$(V^d, V^d)'$  and this is now clear.

 On the other hand by b) of Lemma 2.6,
 $\eta_V\in(V, V)'\otimes(K, K)\ .$ To prove the claim it remains to show 
that the first component of $W_{12}(\eta_V)_{13}{W_{12}}^{-1}$ is in 
$(V, V)'$. Now 
$$1_{KH}\times K^*W_{12}(\eta_V)_{13}{W_{12}}^{-1}1_{KH}\times K=
W1_{KH}\times K^*(\eta_V)_{13}1_{KH}\times K{W}^{-1}$$
$$\subset W(V, V)'\otimes{\mathbb C}W^{-1}
\subset (V, V)'\otimes(W, W)'$$ by the corepresentation relation, 
as claimed.
Let  $\eta_{W}$ be the unitary on $HK$ defined by $(\eta_{W})_{23}=X_W\ .$
A straightforward computation shows that $X_W\in(V\times W\times V,
V\times \iota(W)\times V)\ ,$
thus $\eta_{W}\in(W\times V, \iota(W)\times V)\ ,$ and that 
$W\in\C(V)\mapsto\eta_{W}$ is a natural transformation from $R$ to $R\iota$.
We now check that $(2.3)$ holds.
%
 $$(\eta_{W\times W'})_{234}=X_{W\times W'}=$$
$${{W'}_{43}}^{-1}{{W}_{42}}^{-1}(\eta_{V}^{-1})_{14}{W}_{42}{W'}_{43}W_{12}{W'}_{13}(\eta_{V})_{14}
{{W'}_{13}}^{-1}{W_{12}}^{-1}=$$
$${{W'}_{43}}^{-1}{{W}_{42}}^{-1}({\eta_{V}}^{-1})_{14}{W}_{42}W_{12}
(\eta_{V})_{14}{W'}_{43}(\eta_{W'})_{34}{W_{12}}^{-1}=$$
$${{W'}_{43}}^{-1}(\eta_{W})_{24}{W'}_{43}(\eta_{W'})_{34}\ .$$
Finally, we prove the last statement. Let us assume that
$\varepsilon$ is standard on ${\cal C}(V)_V\ ,$ so  a) of Lemma 2.6 holds.
Hence $\var$ is standard on ${\cal C}(V)$ by Lemma 2.4, completing the
proof.
 \medskip

  We now describe one way of getting  standard braided
symmetries on ${\cal C}(V)$.\medskip

\noindent{\bf 2.11 Proposition} {\sl Let $V$ be a  multiplicative unitary
and $U\in{\cal U}(K)$ such that 
$$\hat V=I\times U\vartheta V\vartheta I\times U^*$$
 is multiplicative, with $\vartheta=\vartheta_{K, K}$. If 
$[\hat{V}_{12}, V_{23}] =0$ and $W\in{\cal C}(V)$, then 
there is a standard braided symmetry $\varepsilon$ on ${\cal C}(V)$ 
defined by:
$$\varepsilon_W:=WU\times IWU^*\times I\vartheta_{H,K}.$$ 
The corresponding natural transformation $\eta$ 
making $V$ into a right regular object is given by   
$$\eta_W: =I\times U\vartheta_{K, H} W\vartheta_{H, K} I\times U^*.$$ }\smallskip

\noindent{\bf Proof.} It is obvious from the form of $\varepsilon$ that we 
have a natural transformation. Hence $\eta$ will be a natural 
transformation, too and a simple computation shows that it makes 
$V$ into a right regular object. Since 
$d')$ of Lemma 2.4 holds, $\var$ is a standard braided symmetry. 
\medskip

If $U\in{\cal U}(K)$ has the properties listed in 
\cite{BS} to make $V$ an irreducible   
multiplicative unitary
 then all the conditions of the above proposition are satisfied. In
particular, ${\cal C}(V)$ has a canonical standard braided symmetry if $V$
comes from a Kac--von Neumann algebra as in \cite{BS} or  is any regular discrete or
compact multiplicative unitary. If $V$ is derived from a locally compact
 group $G$,  the corresponding braided symmetry is that derived from the
usual permutation symmetry on the representation category of $G$ 
interchanging the order of factors in the tensor product of two
representations.\bigskip

\section{Conjugation}

  Our aim is to discuss conjugation in the context of multiplicative 
unitaries and their associated Hopf algebras. Although this aspect was not 
discussed in \cite{BS}, 
relevant related work can be found in a number of 
publications,  and we refer, in particular, 
to the work of Woronowicz 
in the context of compact quantum groups \cite{WTK}.\smallskip    

  However, some of the relevant problems can be seen at the level of the 
representation theory of $C^*$--algebras and von Neumann algebras and it is 
hence wise to discuss them in this simplified setting. We therefore begin 
with $C^*$--categories and $W^*$-categories. If $\cal{T}$ is a $C^*$--category, 
then a {\it conjugation} on $\cal{T}$ is an extension ${\cal T}^a$ of $\cal{T}$ with 
the same objects to include antilinear arrows with the property that any object 
is the source of an antiunitary. To formalize the structure involved, we define a 
{\it semilinear} $C^*$--category to be a $C^*$-category where for each pair 
of objects $\rho$, $\sigma$ in addition to the linear space $(\rho,\sigma)$ 
of ``linear'' arrows there is a second linear space $(\rho,\sigma)_a$ of 
``antilinear'' arrows. The composition of two arrows is antilinear if
and only if precisely one of them is antilinear. Identity arrows are,  
of course, linear and we have 
$$\mu s\circ\lambda r=\mu\lambda s\circ r,$$
$$\mu s\circ\lambda r=\mu\ov\lambda s\circ r,$$ 
according as $s$ is linear or antilinear. The adjoint $r\mapsto r^*$ is 
a contravariant involution leaving objects fixed and being antilinear 
on linear arrows and linear on antilinear arrows. The spaces 
$(\rho,\sigma)$ and $(\rho,\sigma)_a$ are equipped with a norm making 
them into Banach spaces and having the $C^*$--property:
$$\|r\|^2=\|r^*\circ r\|.$$ 
If we forget the antilinear arrows, we get an ordinary $C^*$--category 
and the norm is determined by its values on that subcategory.\smallskip 

   An antiunitary arrow in a semilinear $C^*$--category is an arrow $J$ 
in some $(\rho,\sigma)_a$ such that $J^*\circ J=1_\rho$ and 
$J\circ J^*=1_\sigma$. Two objects $\rho$ and $\ov\rho$ are said 
to be {\it conjugates} if there exists an antiunitary
$J\in(\rho,\ov\rho)_a$. 
Conjugates, if they exist, are defined up to unitary equivalence.\smallskip

  The above definition would seem to be the most natural from 
the categorical point of view. However, if conjugates exist, we may wish to make 
a choice, $\rho\mapsto J_\rho$, of antiunitary 
for each object $\rho$ and then there is an associated antilinear  
$^*$--functor on $\cal{T}$ defined by 
$$\ov T:=J_\sigma\circ T\circ J_\rho^*\in(\ov\rho, \ov\sigma),\quad
T\in(\rho,\sigma).$$ 
It can be extended to ${\cal T}^a$ by defining 
$$\ov{R\circ J_\rho}:=\ov R\circ J_{\ov\rho}$$ 
on antilinear arrows. In addition there is an associated natural unitary 
transformation $d_\rho:(\rho,\ov{\ov\rho})$ defined by
$$d_\rho:=J_{\ov\rho}\circ J_\rho$$ 
and satisfying $\ov{d_\rho}=d_{\ov\rho}$.\smallskip 

   More interestingly, we can also go in the other direction. If we are 
given an antilinear $^*$--functor and a unitary 
natural transformations $d$, as above, we may define 
a semilinear $C^*$--category as follows. For each object $\rho$, 
we introduce an antiunitary arrow $J_\rho\in(\rho,\ov\rho)_a$. A general 
antilinear arrow in $(\rho,\sigma)_a$ can now be written uniquely in the 
form $R\circ J_\rho$, where $R\in(\ov\rho,\sigma)$. Composition with a 
linear arrow $P\in(\pi,\rho)$ is defined by 
$$R\circ J_\rho\circ P:=R\circ\ov P\circ J_\pi.$$ 
Composition with a linear arrow $S\in(\sigma,\tau)$ is defined by 
$$S\circ R\circ J_\rho:=(S\circ R)\circ J_\rho.$$ 
Finally, composition with an antilinear arrow $S\circ J_\sigma$, where 
$S\in(\ov\sigma,\tau)$, 
is defined by 
$$S\circ J_\sigma\circ R\circ J_\rho:=S\circ\ov R\circ d_\rho.$$ 
Routine computations verify that we get a $^*$--category and indeed a 
semilinear $C^*$--category if we define the norms of antilinear arrows in the only 
way compatible with $J_\rho$ being antiunitary, namely by setting 
$$\|R\circ J_\rho\|:=\|R\|.$$
 
  It should be noted that in the above construction of ${\cal T}^a$ if 
$U_\rho\in(\ov\rho,\tilde\rho)$ is a unitary natural transformation 
between two antilinear $^*$--functors then mapping antilinear arrows 
by $R\circ\tilde J_\rho\mapsto R\circ U_\rho\circ J_\rho$ and 
leaving linear arrows invariant is an isomorphism of the constructed 
semilinear tensor $C^*$--categories. Two different choices, $\rho\mapsto J_\rho$ 
and $\rho\mapsto \tilde J_\rho$, within ${\cal T}^a$ lead to a unitary 
natural equivalence $U_\rho:=\tilde J_\rho\circ J_\rho^*$ between 
the associated antilinear $^*$--functors.\smallskip 
    
  We may want our conjugation to have additional properties. The following 
definition would seem 
to describe the best possible situation. We call a {\it strict involutive 
conjugation} an involutive antilinear covariant functor on $\cal{T}$ 
commuting with the adjoint, 
taking an object $\rho$ to $\ov\rho$ and an arrow $T$ to $\ov T$. If we 
now adjoin to the category, as a special case of the above construction, 
an antiunitary $J_\rho$ for each object $\rho$ 
with the property that $J_{\ov\rho}=J_\rho^*$ and
$$J_\sigma T=\ov TJ_\rho,\quad T\in(\rho,\sigma)$$  
then we will have constructed a conjugation on $\cal{T}$. This special case 
corresponds to being able to take $d$ as the identity natural transformation. 
Looked at from the point of view of ${\cal T}^a$, it means that $J_\rho$ can be 
chosen so that $J_{\ov\rho}=J_\rho^*$.
\smallskip 

  To give a simple example: let $\cal{H}$ be a category of Hilbert spaces then 
we get a conjugation on $\cal{H}$ by adding to the arrows all bounded 
antilinear mappings between the respective objects. Such a category will be 
denoted ${\cal H}^a$ and referred to as a category of Hilbert spaces with 
conjugation. Pick an 
orthonormal basis for each Hilbert space $H$ in the category and let
$J_H$ denote the antiunitary involution on $H$ leaving this basis fixed. Then 
define for $T\in(H,K)$, $\ov T:= J_KTJ_H$ and we have a strict involutive 
conjugation on $\cal{H}$ yielding ${\cal H}^a$ as the associated conjugation.\smallskip 

   A second simple example is provided by  a $C^*$--algebra $\cal{A}$ equipped 
with a conjugation $j$, i.e.\ an antilinear involutive $^*$--homomorphism. 
Consider the representation theory of $\cal{A}$ on the objects of $\cal{H}$. 
If $\pi$ is such a representation, we write  $J_\pi:=J_{H_\pi}$ and define 
$$\ov\pi(A):=J_\pi\pi(j(A))J_{\ov\pi},\quad A\in\cal{A},$$ 
$$\ov T=J_\sigma TJ_\rho,\quad T\in(\rho,\sigma).$$ 
In this way, we get a strict involutive conjugation on the $C^*$--category of 
representations of $\cal{A}$ on the objects of $\cal{H}$. The forgetful 
functor into $\cal{H}$ preserves the strict involutive conjugation 
in the obvious sense.\smallskip 

   There is also a simple result going in the other direction. We recall 
\cite{GLR} 
that if $H:\cal{T}\to\cal{H}$ is a $^*$--functor of a $C^*$--category into the 
category of Hilbert spaces then the bounded natural transformations from $H$ 
to $H$ form a von Neumann algebra denoted $(H,H)$ and called the 
commutant of $H$. The evaluation maps $\eta\mapsto\eta_\rho$ are normal 
representations of $(H,H)$. When $\cal{T}$ is a $W^*$--category and $H$ is 
faithful and normal, then $\cal{T}$ can be 
interpreted as a category of normal representations of $(H,H)$.
\smallskip 

\noindent
{\bf 3.1 Lemma} {\sl Let ${\cal T}^a$ be a conjugation on $\cal{T}$ and
$H^a:{\cal T}^a\to{\cal H}^a$ be a $^*$--functor into a category of 
Hilbert spaces with conjugation and $H$ its restriction to $\cal{H}$. 
Given $\eta\in(H,H)$, set 
$$j(\eta){}_\rho:=H^a(J_\rho^*)\eta_{\ov\rho} H^a(J_\rho),$$ 
where $J_\rho$ is an antiunitary from $\rho$ to $\ov\rho$ in ${\cal T}^a$. 
Then $j$ is a conjugation on $(H,H)$.}\smallskip  

\noindent
{\bf Proof.} Given $T\in(\rho,\sigma)$, 
$\ov T:=J_\sigma TJ_\rho^*\in(\ov\rho,\ov\sigma)$ 
and a simple computation shows that $j(\eta)\in(H,H)$. Two different 
choices of $J_\rho$ differ by a unitary in $\cal{T}$. But $\eta$ is a 
natural transformation, so $j$ does not depend on the choice of $J_\rho$ 
and this makes it obvious that $j$ is an involution.\smallskip 

   We next show how, given a conjugation on a $C^*$--algebra, the 
GNS construction provides canonical antiunitary intertwiners between 
conjugate representations.\smallskip 

\noindent
{\bf 3.2 Lemma} {\sl Let $j$ be a conjugation on a $C^*$--algebra $\cal{A}$ 
and let $\phi$ denote a lower semicontinuous densely defined weight on 
$\cal{A}$ and let $\ov\phi:=\phi\circ j$. Let ${\cal N}_\phi$ and 
${\cal N}_{\ov\phi}$ be the associated scalar product spaces mapped by
$\ \ \hat {}\ \ $ into the associated Hilbert spaces,  
$L^2(\cal{A},\phi)$ and $L^2({\cal A},\ov{\phi})$, respectively. 
Then there is a canonical antiunitary operator $J_\phi$ 
from $L^2(\cal{A},\phi)$ to $L^2({\cal A},\ov{\phi})$ defined by 
$$J_\phi\hat N=\widehat{j(N)},\quad N\in{\cal N}_\phi,$$ 
and we have 
$$J_\phi\pi_\phi(A)=\pi_{\ov\phi}\circ j(A)J_\phi,\quad A\in\cal{A}.$$
If $\phi$ extends to a faithful normal weight on $\pi_\phi(\cal{A})''$ 
then 
$$J_\phi S_\phi=S_{\ov\phi}J_\phi,$$
where the operators $S$ are the closed operators derived from the adjoint.}\smallskip

\noindent
{\bf Proof.} $J_\phi$ is uniquely defined as an antiunitary operator since 
$$\phi(N^*N)=\ov\phi(j(N)^*j(N)),\quad N\in{\cal A}.$$ 
Furthermore, the intertwining property holds since
$$j(AN)=j(A)j(N),\quad A\in{\cal A},\,\,N\in{\cal N}_\phi.$$
The final relation follows from 
$$j(N^*)=j(N)^*,\quad N\in{\cal N}_\phi.$$






\section{Conjugation and Tensor Products}
 
   After this general discussion of conjugation   which already illustrates 
the basic problems involved, we turn to conjugation on tensor 
$C^*$--categories, the structures arising in the representation theory of 
Hopf $C^*$--algebras, locally compact quantum groups and multiplicative 
unitaries.\smallskip 

   In a semilinear tensor $C^*$--category, the tensor product is 
defined separately for linear and antilinear arrows. If $R\in(\rho,\sigma)$ and 
$R'\in(\rho',\sigma')$ then $R\times R'\in(\rho\rho',\sigma\sigma')$ and the 
map $R,R'\mapsto R\times R'$ is bilinear. If $R\in(\rho,\sigma)_a$ and 
$R'\in(\rho',\sigma')_a$, then $R\times R'\in(\rho\rho',\sigma'\sigma)_a$ 
and $R,R'\mapsto R\times R'$ is again bilinear.  If $S\in(\sigma,\tau)$ 
and $S'\in(\sigma',\tau')$, then 
$$S'\times S\circ R\times R'=(S\circ R)\times (S'\circ R').$$ 
If $P\in(\pi,\rho)$ and $P'\in(\pi',\rho')$ then 
$$R\times R'\circ P\times P'=(R\circ P)\times (R'\circ P').$$ 
If $S\in(\sigma,\tau)_a$ and $S'\in(\sigma',\tau')_a$ then again 
$$S'\times S\circ R\times R'=(S\circ R)\times(S'\circ R').$$ 
Finally, $(R\times R')^*=R^{'*}\times R^*$ for $R$, $R'$ antilinear.\smallskip 

   A simple example of a semilinear tensor $C^*$--category is got by considering the 
linear and antilinear intertwining operators between a set of unitary 
representations of a group $G$, closed under tensor products, where an 
antilinear intertwining operator $R\in(\rho,\sigma)_a$ is a bounded antilinear 
operator from $H(\rho)$ to $H(\sigma)$, the underlying Hilbert spaces, 
such that 
$$R\rho(g)=\sigma(g)R,\quad g\in G.$$ 
The only point to note is that the tensor product $R\times R'$ is not 
the usual tensor product $R\otimes R'$ of antilinear operators, but is given 
by 
$$R\times R'=\vartheta(\sigma,\sigma')\circ R\otimes R'=R'\otimes R\circ\vartheta(\rho,\rho')$$ 
where $\vartheta$ is the symmetry on the underlying tensor $C^*$--category of 
Hilbert spaces.\smallskip 

  The idea of conjugation in \S 3 nows adapts 
to a tensor $C^*$--category $\cal{T}$. It is an extension ${\cal T}^a$ of 
$\cal{T}$ to a semilinear tensor $C^*$-category where every object is 
the source of an antiunitary. At this point we make contact with the 
notion of conjugation introduced by  Hayashi and Yamagami\cite{HY}. If 
we pick antiunitaries $J_\rho$ for each object $\rho$ and set 
$$c_{\rho,\sigma}:=J_{\rho\sigma}\circ (J_\rho\times J_\sigma)^*,$$ 
we get a natural unitary equivalence from $\ov\sigma\ov\rho$ to 
$\ov{\rho\sigma}$. A computation shows that 
$$c_{\tau,\sigma\rho}\circ c_{\sigma,\rho}\times 1_{\ov\tau}
=c_{\tau\sigma,\rho}\circ 1_{\ov\rho}\times c_{\tau,\sigma},$$  
$$\ov c_{\rho, \sigma}\circ c_{\ov\sigma,\ov\rho}\circ d_\rho\times
d_\sigma= 
d_{\rho\sigma},$$
where $d$ is the natural unitary equivalence of \S 3. Conversely, 
given an antilinear functor $T\mapsto\ov T$ and the natural equivalences 
$c$ and $d$ satisfying the above relations, then the semilinear 
$C^*$--category $\T^a$ constructed in \S 3 can be made into a 
tensor $C^*$--category by using the following definition of the 
tensor product of antilinear arrows: 
$$(R\circ J_\rho)\times (R'\circ J_{\rho'}):=R'\times R\circ c_{\rho,\rho'}^{-1}
\circ J_{\rho\rho'}.$$

  Rather than using semilinear structure, Hayashi and Yamagami define 
a conjugation as an antilinear $^*$--functor equipped with the 
natural transformations $c$ and $d$. \smallskip 
 
  They also introduce the notion of a strict conjugation on a tensor 
$C^*$--category requiring $c$ and $d$ to be identities. In terms of 
antiunitary operators, this obviously corresponds to requiring that 
there is a choice of $J$ such that
$$J_{\rho\sigma}=J_\rho\times J_\sigma\in(\rho\sigma,\ov\sigma\ov\rho)$$ 
for each pair of objects $\rho$ and $\sigma$ and we refer to this case 
as a strict involutive conjugation of tensor $C^*$--categories.\smallskip

  We give an example of a strict tensor $W^*$--category of Hilbert spaces 
with conjugation. Let $\cal{M}$ be a von Neumann algebra equipped with 
a conjugation $j$. Let the objects of the category be the Hilbert spaces 
in $\cal{M}$. If $H$ is such a Hilbert space then its conjugate is $j(H)$. 
If $T\in(H,H')$ then its conjugate is $j(T)\in(j(H),j(H'))$. As this 
conjugation is involutive, we may take the natural unitary equivalence $d$ 
to be the identity. The natural unitary equivalence $c$ from $j(K)j(H)$ 
to $j(HK)$ is defined by 
$$c_{H,K}:=\theta(j(K),j(H)).$$ 
Since $j(c_{K,H})\circ c_{j(K),j(H)}=1_{KH}$, we may construct a 
semilinear tensor $W^*$--category with conjugation, as explained above.
This construction is realized concretely by taking as antiunitary 
arrows $J_H$, the mapping $\psi\mapsto j(\psi)$ for $\psi\in H$, and 
defining 
$$J_H\times J_K:=c_{H,K}^{-1}\circ J_{HK}.$$ 
 
   For a second example, the category of matrices with complex entries is a 
$C^*$--category in a natural way and becomes a strict tensor $C^*$--category 
when the tensor product is defined using lexicographical ordering. 
However, this cannot be made into a strict tensor $C^*$--category with 
a strict involutive conjugation. In fact, labelling the objects by 
the integers in the obvious way, the equation 
$$J_2\times J_3=J_3\times J_2$$ 
cannot be satisfied. On the other hand, our axioms require
$\ov\sigma\ov\rho$ 
rather than $\ov\rho\ov\sigma$ to be the conjugate of $\rho\sigma$. If 
we use the ordinary tensor product of antilinear operators, denoted by 
$\otimes$, then we can satisfy 
$$J_m\otimes J_n=J_{mn}=J_n\otimes J_m,\quad m,n\in{\Bbb N}$$ 
by defining 
$$J_me_i:=e_{m+1-i},$$ 
with respect to the natural orthonormal basis $e_i$, $i=1,2,\dots
m$.\smallskip

  For a third example, we consider a strict tensor $C^*$--category with 
conjugates
\cite{LR} embedded in a strict tensor category of Hilbert
spaces. Let 
$R\in(\iota,\ov \rho\rho)$ and $\ov R\in(\iota,\rho\ov\rho)$ solve the 
conjugate equations for an object $\rho$ of $\cal{T}$, then they also solve 
the corresponding equations in the category of Hilbert spaces and there
is an 
 invertible antilinear operator $T$  from 
$H$ to $\ov H$, the underlying Hilbert space of $\rho$ and $\ov\rho$, 
such that 
$$R(1)=\sum_iTe_i\otimes e_i,\quad \ov R(1)=\sum_ie_i\otimes T^{*-1}e_i,$$ 
where $e_i$ is an orthonormal basis of $H$. We set $\ov T:={T^*}^{-1}$. 
If we pick, for each object $\rho$ of $\cal{T}$, a standard solution 
of the conjugate equations and denote the antilinear operator by $T_\rho$, 
then\cite{RT} there is an 
antilinear functor $S\mapsto \ov S$ commuting with the adjoint defined by 
$$\ov S:=T_\sigma ST_\rho^{-1},\quad S\in(\rho,\sigma),$$
$$f_\rho:=T_\rho^*T_\rho$$ 
is independent of the choice of $T_\rho$, $\rho\mapsto f_\rho$ is a natural 
transformation of the embedding functor to itself and 
$$f_{\rho\sigma}=f_\rho\otimes f_\sigma.$$
One can similarly define an antilinear functor $S\mapsto\tilde S$
associated to the antilinear operators $\ov{T}_\rho:={{T_\rho}^*}^{-1}$
by
$$\tilde{S}:=\ov{T}_\sigma S {\ov{T}_\rho}^{-1}, \quad S\in(\rho,
\sigma).$$

\noindent
{\bf Theorem 4.1} {\sl
\begin{description}
\item{\rm a)} 
Setting $$(\rho, \sigma)_a:=(\ov\rho, \sigma)\circ T_\rho$$
and 
$$T_\rho\times T_\sigma:=T_\sigma\otimes T_\rho\circ\theta(H,K),$$ 
where $H$ and $K$ denote the underlying Hilbert spaces of $\rho$ and $\sigma$, 
respectively, gives an embedded semilinear tensor category.
Set $d_\rho:=T_{\ov\rho}\circ T_\rho$ and $c_{\rho,\sigma}:=T_{\rho\sigma}
\circ(T_\rho\times T_\sigma)^{-1}$.
Then $d$ and $c$ are unitary natural transformations 
satisfying the identities needed to yield a semilinear tensor $C^*$--category 
with conjugation $\T^a$. 

\item{\rm b)} Setting $$(\rho, \sigma)_{\ov a}:=(\ov\rho,
\sigma)\circ\ov T_\rho$$
gives, in a similar way, another semilinear tensor category with
conjugation ${\cal T}^{\ov a}$ with the corresponding properties. 
\end{description}

In
general, $T_\rho\notin(\rho, \ov\rho)_{\ov a}$.
 The category ${\cal T}^a$ (resp. ${\cal T}^{\ov a}$)  may be
identified with the embedded semilinear tensor category after  
the adjoint of antilinear arrows has been redefined so as to make 
the $T_\rho$ (resp. $\ov T_{\rho}$) antiunitary. It is independent of the
embedding of $\cal{T}$ 
into a tensor $C^*$--category of Hilbert spaces.}\smallskip  

\noindent
{\bf Proof.} The observation on the embedded semilinear category is already 
made in \cite{RT}. We get a semilinear tensor category since 
$T_\rho\times T_\sigma$ would be a possible choice of $T_{\rho\sigma}$. 
Since $S\mapsto \ov S$ is defined in terms of $\rho\mapsto T_\rho$, 
$d$ and $c$ are obviously natural transformations satisfying the required 
identities. $d_\rho$ is unitary since $T_\rho^{-1}$ would be a possible 
choice of $T_{\ov\rho}$. Similarly, $c_{\rho,\sigma}$ is unitary since 
$T_\rho\times T_\sigma$ is a possible choice of $T_{\rho\sigma}$. The 
semilinear tensor $C^*$--category determined by this data is 
obviously 
${\cal T}^a$, where the adjoint of antilinear arrows has been changed to 
make the $T_\rho$ antiunitary. Since different choices of the $T_\rho$ 
differ by a unitary, the resulting category is independent of the 
embedding. The statements relative to the category ${\cal
T}^{\ov a}$ can be proved similarly.\smallskip 

   If the semilinear tensor $C^*$--category with conjugation constructed 
above can be embedded in the semilinear tensor category of Hilbert spaces, 
then the $T_\rho$ are antiunitary for this embedding and the intrinsic 
dimensions coincide with the dimensions of the underlying Hilbert spaces.\smallskip

  We can learn more from the above construction of a conjugation. 
The natural transformations $c$ and $d$ have here been defined in terms of 
the invertible antilinear operators $T_\rho$ which in turn were defined using 
standard solutions $R_\rho$ and $\overline R_\rho$ of the 
conjugate equations. Expressing $c$ and $d$ in terms of the $R_\rho$ and 
$\overline R\rho$, we find 
$$c_{\rho,\sigma}=1_{\overline{\rho\sigma}}\times(\overline R^*_\rho\circ 
1_\rho\times\overline R^*_\sigma\times 1_{\overline\rho})\circ 
R_{\rho\sigma}\times 1_{\overline\sigma\overline\rho},$$
$$d_\rho=1_{\overline{\overline\rho}}\times R^*_\rho\circ 
R_{\overline\rho}\times 1_\rho.$$ 
These expressions no longer make reference to an ambient 
Hilbert space. Defining the conjugate linear $*$--functor 
by 
$$\overline S\times 1_\rho\circ R_\rho=1_{\overline\sigma}\times S^*\circ R_\sigma,
\quad S\in(\rho,\sigma),$$ 
$c$ and $d$ become natural transformations. Furthermore, 
after a somewhat lengthy calculation, the identities between 
$c$ and $d$ can be verified, leading to the following 
result.\smallskip 

\noindent
{\bf 4.2 Theorem} {\sl Any strict tensor $C^*$--category with 
conjugates admits a canonical 
conjugation defined, as above, in terms of standard solutions of 
the conjugate equations.}\smallskip 

\noindent
{\bf Proof.} The only point still to be checked is that 
the conjugation does not depend on the choice of standard 
solutions of the conjugate equations. However, a second 
choice $\rho\mapsto\tilde R_\rho$ is related to the first by 
$\tilde R_\rho=U_\rho\times 1_\rho\circ R_\rho$, where 
$U_\rho\in(\overline\rho,\tilde\rho)$ is unitary. We then 
have $\tilde d_\rho=U_{\overline\rho}\circ\overline{U_\rho}$ 
and $\tilde c_{\rho,\sigma}=U_{\rho\sigma}\circ c_{\rho,\sigma}
\circ(U_\sigma\times U_\rho)^*$. As we have seen this leads to 
the same conjugation.\smallskip 

  The reader's attention is drawn to a result of Yamagami's, 
Theorem 3.6 of \cite{Yam}, where he achieves more, at the cost 
of passing to an equivalent tensor $C^*$--category in the course 
of the proof. We also remark that 
$$c_{\rho,\sigma}=1_{\overline{\rho\sigma}}\times(R^*_\sigma\circ 
1_{\overline\sigma}\times R^*_\rho\times 1_\sigma)\circ 
1_{\overline{\sigma}\overline\rho}\times\overline R_{\rho\sigma},$$ 
$$d_\rho=\overline R^*_\rho\times 1_{\overline{\overline\rho}}
\circ 1_\rho\times\overline R_{\overline\rho}.$$

   We recall from the beginning of the previous section that, in the presence of 
a conjugation $j$ on a $C^*$--algebra $\cal{A}$, the representation theory 
relative to a category of Hilbert spaces with a strict involutive conjugation
has a strict involutive conjugation given by 
$$\bar\pi(A):=J_\pi\pi(j(A))J_{\bar\pi},\quad A\in\cal{A},$$ 
$$\bar T=J_\sigma TJ_\rho,\quad T\in(\rho,\sigma),$$ 
where $J_\pi=J_{H_\pi}$ and $H_\pi$ is the Hilbert space of $\pi$.
If $\cal{A}$ is a Hopf $C^*$--algebra with coproduct $\delta$ satisfying 
$$\delta\circ j=j\otimes j\circ\theta\circ\delta,$$ 
and we consider representations relative to a strict tensor $W^*$--category 
of Hilbert spaces with a conjugation, 
then $J_{\pi\times\rho}=J_{\rho}\otimes J_\pi\vartheta(H_\pi, H_{\rho})$ 
defines $\ov{\rho}\times \ov \pi$ as a 
conjugate for $\pi\times \rho$ and we get a conjugation 
on the tensor $W^*$--category of representations of $\cal{A}$. 
If the underlying category of Hilbert spaces has a strict involutive 
conjugation then the same is true for the category of representations 
of $\cal{A}$.\smallskip 

  We now come to other cases where conjugates can be 
defined in terms of antiunitary arrows but where we 
need to make a simple extension of our formalism. 
Instead of starting with a strict tensor $C^*$--category,
we need to adjoin antilinear 2-arrows to a
$2-C^*$--category. A formal definition of $2-C^*$--category 
can be found in \cite{LR} but the examples given below
should be self--explanatory. We consider a set of von
Neumann algebras. These form the 0--arrows. The
bimodules (correspondences)  on this set form the
$1$--arrows, whilst the bimodule  homomorphisms form the
$2$--arrows. Compositions are defined  in the obvious
manner. What we get is not a $2$--$C^*$--category but
what might be called a bi--$C^*$-category, because the 
composition of $1$--arrows is defined only up to
equivalence. Now there is no problem in adjoining  antilinear
$2$--arrows in a natural way because there is  a natural
notion of an antilinear bimodule homomorphism.  An
antilinear bimodule homomorphism from an
$\cal{M}$--$\cal{N}$--bimodule to an
$\cal{N}$--$\cal{M}$--bimodule is simply a bounded
antilinear map $A$ between the underlying Hilbert spaces 
such that 
$$A(M\cdot\psi\cdot N)= N^*\cdot(A\psi)\cdot M^*,\,
M\in{\cal{M}},\, N\in\cal{N}.$$ 
Adding these antilinear $2$--arrows, we get a semilinear 
bi--$C^*$--category, where every $1$--arrow is the source 
of an antiunitary $2$--arrow. In fact, conjugating an 
$\cal{M}$--$\cal{N}$--bimodule with an antiunitary operator
yields an $\cal{N}$--$\cal{M}$--bimodule, a conjugate 
bimodule unique up to equivalence.\smallskip 
  
  In the second example, we deal with morphisms of 
von Neumann algebras and whilst it is not immediately 
evident that we can define antilinear intertwining 
operators between such morphisms, the close links 
between morphisms and bimodules suggest that it 
must be possible. Furthermore, there is a definition of 
conjugation for such morphisms going back to Longo. 
These considerations lead us to consider a separable Hilbert
space and a set of  von Neumann algebras represented
standardly on that  Hilbert space. We denote by
$J_{\cal{M}}$ the corresponding  modular conjugation of the
von Neumann algebra $\cal{M}$ and let $j_{\cal{M}}$ denote
Ad$J_{\cal{M}}$. Then if $\rho:\cal{M}\to\cal{N}$ and 
$\sigma:\cal{N}\to\cal{M}$. Then we write
$A\in(\rho,\sigma)_a$ if $A$ is a (bounded) antilinear 
operator on our Hilbert space such that 
$$A\rho(M)=j_{\cal{M}}(M)A,\quad M\in\cal{M},$$ 
$$Aj_{\cal{N}}(N)=\sigma(N)A,\quad N\in\cal{M}.$$ 
As these conditions may look surprising, it is perhaps 
worth observing that if $\tau:\cal{M}\to\cal{N}$ then 
the condition for a (bounded) linear operator $T$ to be 
in $(\rho,\tau)$ is that 
$$T\rho(M)=\tau(M)T,\quad M\in\cal{M},$$ 
$$Tj_{\cal{N}}(N)=j_{\cal{N}}(N)T,\quad N\in\cal{N}.$$ 
Composition of these $2$--arrows can be defined in the 
obvious fashion. The tensor product for linear $2$--arrows 
is well known. For antilinear arrows, we proceed as 
follows: if $A$ is as above and $A'\in(\rho',\sigma')_a$, 
where $\rho':\cal{P}\to\cal{M}$ and
$\sigma':\cal{M}\to\cal{P}$, then 
$$A\times
A':=A'J_{\cal{M}}A\in(\rho\rho',\sigma'\sigma)_a.$$
It is easy to verify that we get a semilinear
$2$--$C^*$--category in this way. It is also easy to 
recover Longo's result on the existence of conjugates. 
Given $\rho$ as above, we want to show that $\rho$ is the
source of an antiunitary arrow $A$. The first of the
equations that $A$ has to satisfy can be solved by taking 
$A:=J_{\cal{M}}U^*$, where $U$ is a unitary implementing
$\rho$. We now set 
$$\overline\rho(N):=Aj_{\cal{N}}(N)A^*,$$ 
and can check that
$\overline\rho:\cal{N}\to\cal{N}$.\smallskip 

  We now adapt the above formalism to the case of $C^*$--algebras 
by replacing the above von Neumann algebras by weakly dense unital 
$C^*$--algebras $\A$, $\B$, $\C,\dots$. We write $J_\A$ for the 
modular conjugation of the weak closure of $\A$. 
The difference is that we now no longer have the analogue of 
Longo's result on the existence of conjugates. We therefore 
consider a $2$--$C^*$--category $\T$ whose $0$--arrows are this set of 
$C^*$--algebras whose $1$--arrows are morphisms between these 
$C^*$--algebras and whose $2$--arrows are intertwiners between these 
morphisms. We suppose we may pick for each $1$--arrow $\rho:\A\to\B$ 
an antiunitary operator $J_\rho$ in such a way that 
\begin{description}
\item{a)} $J^*_\rho J_\A$ induces $\rho$, 
\item{b)} $J_\rho J_\B$ induces a $1$--arrow 
$\overline\rho:\B\to\A$.
\end{description} 
It follows that a $1$--arrow $\rho$, being unitarily implemented, 
extends to a morphism $\hat\rho$ between the weak closures and is 
a $1$--arrows in the semilinear $2$--$C^*$--category, $\C_a$ say, 
of morphisms and intertwiners between the weak closures. It follows 
easily from a) and b) that $J_\rho\in(\hat\rho,\hat{\overline\rho})$. 
We give conditions allowing us to construct a conjugation on $\T$ 
as a semilinear $2$--$C^*$--subcategory of $\C_a$. The antiunitaries 
$J_\rho$ are not unique and two choices $J_\rho$ and $\tilde J_\rho$ 
are said to be equivalent if $J^*_\rho\tilde J_\rho\in(\rho,\rho)$ 
and $J_\rho\tilde J_\rho^*\in(\overline\rho,\overline\rho)$. The 
equivalence class of $J_\rho$ is denoted by $[J_\rho]$. We now require  
\begin{description} 
\item{c)} $J^*_\rho\in[J_{\overline\rho}]$, 
\item{d)} $J_\rho J_\B J_{\rho'}\in[J_{\rho'\rho}]$,
\item{e)} If $R\in(\pi,\rho)$ in $\T$, where $\pi$ and $\rho$ are 
$1$--arrows from $\A$ to $\B$, then $J_\rho R J_\lambda^*\in\A$. 
\end{description} 
We now claim that we get a conjugation $\T_a$ on $\T$ by taking 
$(\rho,\sigma)_a$ to be the set of antilinear operators on the 
underlying Hilbert space such that $J_\sigma S\in(\rho,\overline\sigma)$ 
and $SJ_\rho^*\in(\overline\rho,\sigma)$. Note that this definition is 
independent of the choice of representatives in the equivalence 
classes and that $(\rho,\sigma)_a\subset(\hat\rho,\hat\sigma)$. 
It is easily checked that the image of $\T_a$ in $\C_a$ is closed under 
$\circ$--composition and adjoints. We note that by c), 
$J_\rho\in(\rho,\overline\rho)$. It therefore remains to show that 
the image of $\T_a$ is closed under tensor products. Let $S\in(\rho,\sigma)_a$ 
and $S'\in(\rho',\sigma')_a$, then working in $\C_a$, we have 
$$J_{\sigma}\times J_{\sigma'}\circ S'\times S=(J_{\sigma'}\circ S') 
  \times(J_{\sigma}\circ S)\in(\rho'\rho,\overline{\sigma'}\overline{\sigma}),$$ 
$$S'\times S\circ(J_{\rho'}\times J_{\rho})^*=(S\circ J^*_{\rho})\times
(S'\circ J_{\rho"}^*)\in(\overline{\rho}\overline{\rho'},\sigma\sigma'),$$ 
so in virtue of d), $S'\times S\in(\rho'\rho,\sigma\sigma')$. Thus  
$\T_a$ is a $2$--$C^*$--subcategory of $\C_a$ and $\T_a$ is a 
conjugation on $\T$.\smallskip

   We now return to the concept of regular object from \S 2 and show that 
the conjugate of a left regular object is a right regular object.\smallskip 

\noindent{\bf 4.3 Proposition} {\sl We consider a strict tensor $W^*$--category 
$\cal{T}$ equipped with a faithful idempotent tensor $^*$--functor $\iota$ 
onto a tensor $W^*$--subcategory of Hilbert spaces. We suppose that 
$\cal{T}$ admits a conjugation $J$ such that for each pair $W$, $W'$ of 
objects of $\cal{T}$ 
$$J_{\iota(W\times W')}\circ J_{\iota(W')}^{-1}\times J_{\iota(W)}^{-1}=
\iota(J_{W\times W'}\circ J_{W'}^{-1}\times J_W^{-1}).$$ 
Then if $V$ is a left regular object of $\cal{T}$, its conjugate $\ov V$ 
is a right regular object.}\smallskip

\noindent{\bf Proof.} Let $R^{\ov V}$ denote the functor of tensoring on 
the right by $\ov V$ and let $\xi\in(L, L\iota)$ denote the
natural unitary transformation making $V$ into a left regular
object. We show that there is  a natural unitary transformation in
$(R^{\ov V}, R^{\ov V}\iota)$.   
We have, for every object $W$ of $V$, an antiunitary arrow $J_W$  defining 
a conjugate $\ov W$ of $W$. 
There is a unitary ${\eta^{\ov V}}_W\in(W\times\ov V,
\iota(W)\times\ov V)$ defined explicitly by
$${\eta^{\ov V}}_W=J_V\times{J_{\iota(W)}}^{-1}\circ\xi_{\ov
W}\circ J_{W}\times{J_V}^{-1}.$$  
$\eta$ is a natural transformation. In fact, if $T\in (W,W')$ then 
$$\eta^{\ov V}_{W'}\circ T\times 1_{\ov V}=J_V\times J_{\iota(W')}^{-1}
\circ \xi_{\ov W'}\circ 1_V\times\ov T\circ J_W\times J_V^{-1}$$ 
$$=J_V\times J_{\iota(W')}^{-1}\circ 1_V\times\iota(\ov T)\circ \xi_{\ov
W}
\circ J_W\times J_V^{-1}=\iota(T)\times 1_{\ov V}\circ\eta^{\ov V}_W.$$  
 We have to show that $\eta^{\ov V}$ satisfies the coherence relation 
$${\eta^{\ov V}}_{W\times W'}=({{\eta^{\ov V}}_{W}})_{13}({{\eta^{\ov V}}_{W'}})_{23}.$$ 
Now $J_W\times J_{W'}\circ J_{W\times W'}^{-1}\in(\overline{W\times W'},
\ov W'\times\ov W)$. Thus by the naturality of $\xi$, we have 
$$\xi_{\ov W'\times\ov W}\circ 1_V\times(J_W\times J_{W'}\circ 
J_{W\times W'}^{-1})=1_V\times\iota(J_W\times J_{W'}\circ J_{W\times W'}^{-1})
\circ\xi_{\overline{W\times W'}}.$$ 
Hence, we may express $\eta^{\ov V}_{W\times W'}$ using $\xi_{\ov
W'\times\ov W}$ 
rather than $\xi_{\overline{W\times W'}}$ and, after simplifying, this yields 
$$\eta^{\ov V}_{W\times W'}=J_V\times J_{\iota(W')}^{-1}\times
J_{\iota(W)}^{-1}
\circ \xi_{\ov W'\times\ov W}\circ J_W\times J_{W'}\times J_V^{-1}.$$  
Finally, writing 
$$\xi_{\ov W'\times\ov W}=(\xi_{\ov W})_{13}\circ J_W\times
J_{\iota(W')}
\times J_V^{-1}\circ J_V\times J_{\iota(W')}^{-1}\times J_W^{-1}\circ
\xi_{\ov W'}
\times 1_{\ov W},$$ 
we get the required result.\smallskip

\section{Conjugation for locally compact quantum groups}

   Conjugation for the representation theory of a locally compact group 
provides motivation for the generalization which follows but is too 
well known to merit discussion. Therefore we turn to consider a
locally compact quantum group, $(\A, \delta, \phi, \psi, R, \tau)$, a concept
which, after initial work of Masuda, Nakagami and Woronowicz,  has perhaps
now received 
its final definition and denomination 
(cf.\  \cite{KvD}, \cite{KV}). Here  $\delta: \A\to
M(\A\otimes \A)$ is a coassociative coproduct such that
$\delta(\A)\A\otimes{\Bbb
C}=\delta(\A){\Bbb C}\otimes\A=\A\otimes\A$, $\phi$ ($\psi$) is called left (right) Haar
measure and 
is a faithful, lower semicontinuous, densely defined, left (right) invariant  KMS weight in
the 
$C^*$--algebraic sense. $R$ is an involutive $^*$--antiautomorphism and 
$\tau$ is a pointwise norm continuous one--parameter automorphism group of 
$\A$ commuting with $R$. $R$ and $\tau$ together should implement the
coinverse, in the
sense that the coinverse of $\A$ should be $\kappa:=R\tau_{-i/2}$, where 
$\tau_{-i/2}$ is the analytic generator of  $\tau$. We focus attention, 
not on the antiautomorphism $R$ but on the conjugation $j:=R\circ *$. 
This defines for us the conjugation on the representation theory and 
makes it clear that the problems of defining a conjugation have been 
avoided by a judicious choice of axioms. On the other hand it must 
be stressed that Woronowicz \cite{W} in effect overcame these difficulties 
in the case of compact matrix pseudogroups by proving that they are 
locally compact quantum groups in a natural way. \smallskip 

  We refer to \cite{KV} for the general definition, recalling here,
instead, 
the properties we need. First, we explain left invariance.  
Let $\M_\phi$ denote as usual the dense $^*$--algebra 
linearly spanned by the $a\in\A^+$ such that $\phi(a)<\infty$ 
and $\N_\phi$ the associated left ideal. 
One has:
$$\phi(\omega\otimes\iota(\delta(a))=\phi(a)\omega(I),\quad \omega\in{\A^*}_+, a\in{\M_\phi}^+.$$
The following property,  referred to as {\it strong left invariance}, is shown 
to hold 
in a locally compact quantum group (\cite{KV}): for all $a,b\in\N_\phi$,
$x:=\overline{\iota\otimes\phi}(\delta(a^*)I\otimes b)$ lies in  the domain of
$\kappa$ and 
$$\kappa(x)=\overline{\iota\otimes\phi}(I\otimes a^*\delta(b)).$$
 $\M_{\iota\otimes\phi}$ 
and $\M_{\overline{\iota\otimes\phi}}$ will denote the 
domains of $\iota\otimes\phi$ and its natural extension to
the multiplier algebra $M(\A\otimes\A)$, cf. 
\cite{K}. In particular we have   
$$\phi(a^*\omega\otimes\iota(\delta(b)))=\phi(\omega\kappa\otimes\iota(\delta(a^*))b
),\quad
a,b\in\N_\phi,$$
with  $\omega$ in the set
   of those elements $\omega\in\A^*$ for which $\omega\kappa\in\A^*$.
The above expression can also be written as
$$\omega\otimes\phi(I\otimes
a^*\delta(b))=\omega\kappa\otimes\phi(\delta(a^*)I\otimes b).$$

Then the following relations also hold.
$$\phi\circ\tau_t=\nu^{-t}\phi, \quad\text{for some } \nu\in\Bbb{R},$$
$$\delta\circ j=j\otimes j\circ\theta\circ\delta$$
$$\tau_t\otimes\tau_t\circ\delta=\delta\circ\tau_t,$$
$\nu$ is referred to as {\it the scaling constant}.

Right invariance can be  derived formally from left invariance choosing
$\phi_r=\phi R$,
 using the antimultiplicativity of $R$ and the relation
$\delta R =\theta R\otimes R\delta$.

  We now turn to multiplicative unitaries starting with a construction of  
\cite{BS}. Let $V$ be a multiplicative unitary 
on $K^2$, and define an associative product on $(K, K)_*$ by
$\omega*\omega'(T)=
\omega\otimes\omega'(V^{-1}I\otimes TV)$,
$T\in(K, K)$ making $(K, K)_*$ into a Banach algebra.\medskip

\noindent{\bf 5.1 Proposition} (\cite{BS}) {\sl For any corepresentation $W$ of $V$,
$\omega\in (K, K)_*\mapsto \omega\otimes i(W)\in(H_W, H_W)$ is an algebra
homomorphism. If $V$ is a regular multiplicative unitary, the closure
$\A(W)$ of its image is a $C^*$--algebra.}\medskip

We endow $(K, K)_*$ with the maximal $C^*$--seminorm determined
by all (unitary) corepresentations of $V$ and denote
the completion of  $(K, K)_*$  in this seminorm by  $\A_{\max}(V)$. 
 We denote by $\pi_W$
the $^*$--representation of $\A_{\max}(V)$ obtained
extending the homomorphism defined in the above proposition.\smallskip

For an operator $X\in(KH,KH)$, we denote 
the norm 
closure of $\{\omega\otimes i(X):\omega\in (K,K)_*\}$
by $\A(X)$  and 
the 
closure of $\{i\otimes \omega(X):\omega\in (H,H)_*\}$ by
$\hat\A(X)$ where $i$ 
denotes the appropriate identity map.
Note that $\omega\in(K, K)_*$ has zero seminorm if and only if $\omega$
annihilates every  $\hat\A(W)$.
Let us assume that $V$ is a regular multiplicative unitary, or, more generally,
that $\hat\A(V)$ is a $C^*$--algebra. Since $W$ is a
corepresentation of $V$, 
$\hat\A(W)$ is contained in $\hat \M(V)$, the von Neumann algebra
generated by $\hat\A(V)$. Therefore $\omega$ will annihilate  
$\hat\A(W)$ if it annihilates $\hat \M(V)$.
Hence $\A_{\max}(V)$ can also be defined as the 
completion of $\hat \M(V)_*$, the predual of $\hat \M(V)$, in the $C^*$--norm
$\|\omega\|=\sup_{W\in{\cal C}(V)}\|\omega\otimes i(W)\|$. \medskip

\noindent{\bf 5.2 Theorem} (\cite{BS}) {\sl If $V$ is a regular multiplicative unitary,
 the map $W\to\pi_W$ defines a faithful 
 $^*$--functor from the $C^*$--category
 ${\cal C}(V)$ of unitary corepresentations of $V$
onto  the $C^*$--category of nondegenerate $^*$--representations  of
$\A_{\max}(V)$.}\medskip

It is also shown in \cite{BS} that $\A_{\max}(V)$ is equipped with a natural
coassociative coproduct $\delta$. One can
easily check, in analogy with the group case, that the map $\rho,\sigma\in
$Rep$(\A_{\max}(V))\to\rho\otimes\sigma\circ\delta\in $Rep$(\A_{\max}(V))$ makes
Rep$(\A_{\max}(V))$ into a tensor $C^*$--category in such a way that $\pi$
is a tensor functor.

   Recalling our aim of defining a conjugate representation
$\ov W$ of a given unitary representation $W$ of $V$, Theorem 5.2 tells 
us that it suffices to determine the associated $^*$--representation 
$\pi_{\ov W}$ of ${\cal A}_{\max}(V)$ and we know from the discussion in
Section 4 that we should look for a suitable conjugation on 
Rep$({\cal A}_{\max}(V))$ or some related $C^*$--algebra. 

If $\A_{\max}(V)$ is equipped with a (densely defined, unbounded) coinverse
$\kappa$, this would be the natural starting point and in view of 
the duality between $\hat\M(V)$ 
and $\hat\M(V)_*$, the coinverse $\kappa'$ and the $^*$--involution of
 $\A_{\max}(V)$ are indeed related by:
$$^*\circ\kappa'(\omega)=\ov{\omega},\quad \omega\in\hat\M(V)_*.$$
There are difficulties involved as 
we shall, in general, be dealing with an antilinear involution that
does not commute with the adjoint and is only densely defined.\smallskip 

   However, it suggests giving a definition in terms of an, in general, unbounded 
antilinear operator. We say that a unitary  $\ov W$ on a 
Hilbert space of the form $KH_{\ov W}$ 
 is a  conjugate of $W$  if there is a densely defined 
closed antilinear invertible operator $T: H_W\to H_{\ov W}$ such that 
$$\omega\otimes i(\ov W)T<T\ov\omega\otimes i(W),\quad \omega\in(K,K)_*.$$

\noindent{\bf 5.3 Lemma} {\sl A conjugate $\ov W$ of a corepresentation 
$W$ is itself necessarily 
a corepresentation.}\smallskip

\noindent{\bf Proof.} 
Let $T:H_W\to H_{\ov W}$ be densely defined closed antilinear 
and invertible, $W$ a corepresentation of $V$ and $\ov W$ 
a unitary on $KH_{\ov W}$ such that 
$$\omega\otimes i(\ov W)T<T\ov\omega\otimes i(W),\quad \omega\in (K,K)_*.$$  
Pick $\phi\in \D_T$ and $\omega_1,\omega_2\in (K,K)_*$ then 
$$\omega_1\otimes\omega_2\otimes i(V_{12}\ov W_{13}\ov
W_{23})T\phi=
T\ov\omega_1\otimes\ov\omega_2\otimes i(W_{13}W_{23}{V_{12}}^*)\phi=$$
$$T\ov\omega_1\otimes\ov\omega_2\otimes i({V_{12}}^*W_{23})\phi=
\omega_1\otimes\omega_2\otimes i(\ov W_{23}V_{12})T\phi.$$ 
Since $T$ is invertible
$$\omega_1\otimes\omega_2\otimes i(V_{12}\ov W_{13}\ov W_{23})=
\omega_1\otimes\omega_2\otimes i(\ov W_{23}V_{12})$$
and since $\omega_1,\omega_2\in(K,K)_*$ are arbitrary, 
$$V_{12}\ov W_{13}\ov W_{23}=\ov W_{23}V_{12}.$$ 
Thus $\ov W$ is a corepresentation of $V$. \medskip

It suffices to verify the intertwining relation
for $\omega$ in a total
set of $(K, K)_*$, e.g.\ 
the set $\{\omega_{\xi, \eta}\}$.
Hence, $\ov W$ is a conjugate of $W$ if and only if for $\eta'\in\D_T$, 
$\xi'\in\D_{T^*}$, $\xi,$ $\eta\in K$,
$$(\xi\otimes\xi', \ov{W}\eta\otimes T\eta')=(W\xi\otimes\eta', \eta\otimes
T^*\xi').\eqno(5.1)$$ 
In practice, this will be checked for $\eta'$ and $\xi'$
in a core of $T$ and $T^*$, respectively. In analogy with the classical 
situation this equation can be interpreted from the standpoint of the Banach 
algebra $\hat{\cal{A}}(V)$ associated with the multiplicative unitary as in 
\cite{BS}: it asserts that the adjoint of the ``matrix coefficient'' 
$i\otimes\omega_{\xi',\eta'}(W)$ is given by 
$i\otimes\omega_{T^{*-1}\xi',T\eta'}(\ov W)$. 
In particular, if the regular representation is 
selfconjugate, $\hat{\cal{A}}(V)$ is $^*$--invariant, hence a $C^*$--algebra. 
We shall see later that the left regular representation of a locally compact 
quantum group is selfconjugate
in this sense. Other examples arise from compact quantum
groups \cite{W}, 
 Hopf--von Neumann algebras \cite{ES} and Kac systems \cite{BS}
as pointed out in \cite{P}. It raises the question of whether
the existence of conjugate representations might prove a
more effective postulate than regularity in developing
the theory of multiplicative unitaries.\smallskip

On the other hand, the above definition of conjugate,
with its reliance on unbounded operators with unspecified
domains, is difficult to work with. It is clear that if
$\overline W$ is a conjugate of $W$ then $W$ is a conjugate
of $\overline W$ since we may use $T^{-1}$ in place of $T$. But
it is not even clear whether a conjugate
is unique up to unitary equivalence. Nevertheless, as we shall
see in the sequel, we can, in special cases, relate this notion
of conjugate to the other notions used in this paper.\smallskip

  We begin by considering the category of finite dimensional
representations of a compact quantum group and add antilinear
operators to get a semilinear category. $T\in(W,\tilde{W})_a$ if $T$ is
an antilinear operator with
$$\omega\otimes i(\tilde W)T=T\overline\omega\otimes i(W),\quad
\omega\in(K,K)_*.$$
Defining tensor products by
$$T\times T':=\theta(\tilde H,\tilde H')\circ T\otimes T',$$
where $T\in(W,\tilde W)_a$ and $T'\in(W',\tilde W')_a$, $H$ and $H'$
are the underlying Hilbert spaces of $W$ and $W'$ and $\tilde H$ and
$\tilde{H}'$ are the underlying Hilbert spaces of $\tilde W$ 
and $\tilde W'$. Adding these
antilinear intertwiners, we get a semilinear tensor
category of bounded intertwiners that is not in general self-adjoint.
We have already met this phenomenon in Theorem 4.1 and can make this
more precise using the antilinear operators $T$ and $\ov T:=T^{*-1}$
associated
with solutions $R$, $\ov R$ of the conjugate equations for $W$ and $\ov
W$
as discussed in conjunction with that theorem. We then have
$$W_{13}\circ 1_K\times R=\ov W_{12}^{-1}\circ 1_K\times R,$$
$$\ov W_{13}\circ 1_K\times\ov R=W_{12}^{-1}\circ 1_K\times\ov R,$$
where $K$ denotes the Hilbert space of the regular corepresentation.
Writing these equations in terms of $T$ instead, we get
$$(\xi\otimes\xi'',W\eta\otimes T^*\xi')=
(\ov W\xi\otimes\xi',\eta\otimes T\xi''),\eqno(5.2)$$
$$(\xi\otimes\xi',\ov W\eta\otimes T^{*-1}\xi'')=
(W\xi\otimes\xi'',\eta\otimes T^{-1}\xi').\eqno(5.3)$$
These equations imply that $T^*\in(\ov W,W)_a$ and $T^{*-1}\in(W,\ov
W)_a$. Thus $\ov W$ is a conjugate of $W$ in the sense of Lemma 5.3, too.
However $T$ is not an intertwiner, 
here reflecting the fact that the antilinear involution $^*\circ\kappa'$ 
on $\A_{\max}(V)$ does not 
commute with the adjoint.\smallskip

  Up till now, we have not explicitly exhibited interesting examples of 
infinite dimensional conjugate corepresentations. It is not surprising that 
this can be done for multiplicative 
unitaries derived from the regular representations of locally compact 
quantum groups.\smallskip

In fact,
the map 
$$a\otimes b\to\delta(b)a\otimes I\text{\quad for\ } a, b\in\N_\phi$$
defines a bounded linear operator $U$ on 
$L^2(\A, \phi)\otimes L^2(\A, \phi)$.
The left invariance of $\phi$
implies that $U$ is isometric:
$$(Ua\otimes b, Uc\otimes d)=
\phi\otimes\phi(a^*\otimes 
I\delta(b^*d)c\otimes
I)=$$
$$\phi(a^*\iota\otimes\phi(\delta(b^*d))c)=
\phi(a^*\phi(b^*d)c)=$$
$$(a\otimes b, c\otimes d).$$
and therefore $U$ is
unitary, since its range is dense (recall that in a locally compact quantum
group
$\delta(\A)\A\otimes I$ 
and $\delta(\A)I\otimes\A$ are assumed to be 
dense in
$\A\otimes\A$).
We next compute the Hilbert space adjoint 
$V:={U}^*$
(which will be  the
standard multiplicative unitary
associated to $\phi$). For $a$, $b$, $c$, $d\in\N_\phi$,
$$(Va\otimes b, c\otimes d)=(a\otimes b, Uc\otimes d)=(a\otimes b,
\delta(d)c\otimes I)=$$
$$\phi(a^*\iota\otimes\phi(I\otimes b^*\delta(d))c)=\phi(a^*\kappa(\iota\otimes\phi(\delta(b^*)I\otimes d))c)$$
by strong left invariance of $\phi$.
Note that both in the proof of being
an isometry and 
in the computation of the adjoint of $U$
 only the left invariance of the second factor of
 $L^2(\A, \phi)\otimes L^2(\A, \phi)$ has been used.

Using a right invariant measure $\phi_r$ (e.g. $\phi_r=\phi\circ R$)
the map $a\otimes b\mapsto\delta(a)I\otimes b$ defines another multiplicative
unitary operator on $L^2(\A, \phi_r)\otimes L^2(\A, \phi_r)$. But we want 
the {\it right regular corepresentation}, $V_r$, instead, a unitary operating 
on  $L^2(\A, \phi)\otimes L^2(\A, \phi_r)$ and this is derived from the map
$a\otimes b\mapsto\vartheta\circ\delta(b)a\otimes 1$, where $\vartheta$ permutes 
the factors in the tensor product.\smallskip

   When we have a locally compact quantum group, we can form a two--parameter
(pointwise norm continuous)
group $\omega_{s,t}=\nu^{\frac{s}{2}}\tau_s\sigma_t$ generated by the modular group $\sigma$
of the left Haar measure $\phi$, and the scaling group $\tau$ (which commutes with $\sigma$). 
It is well known that the set $I(\omega)$ of norm entire elements for $\omega$
is a dense $^*$--subalgebra of  $\A$, 
  stable under
all the maps $\omega_z$, $z\in{\Bbb C}^2$.
Thus $I(\omega)$ is a natural common core for the analytic generators
of both $\sigma$ and $\tau$
\cite{CZ}.
However, we  regard $I(\omega)$ as a subspace of $L^2(\A, \phi)$
so that $\omega$
becomes a unitary
group on  $L^2(\A, \phi)$. We then look for   
  a canonical subspace of $I(\omega)$ dense in   $L^2(\A, \phi)$
which is at the same time a 
common core for  the generators of the  unitary group.
\medskip

\noindent{\bf 5.4 Lemma} {\sl  Let $\phi$ be a lower semicontinuous densely defined KMS weight 
on a $C^*$--algebra $\A$ and let $\omega: {\Bbb R}^n\to\Aut(\A)$ be a pointwise norm continuous 
$\phi$--invariant automorphism group of $\A$ containing the modular group of $\phi$
as a  coordinate subgroup. Then 
${\cal B}_{\phi,\omega}:=\cap_{{z}\in{\Bbb
C}}\omega_{z}\big((\N_\phi\cap{\N_\phi}^*)\cap I(\omega)\big)$, 
is dense in $L^2(\A, \phi)$ and $\omega$ acts on $L^2(\A, \phi)$ as a 
strongly continuous unitary group,  denoted $\Omega$. 
 $\B_{\phi, \omega}$ is a common core
for the positive selfadjoint operators $\Delta_1,\dots, \Delta_n$ 
on $L^2(\A,\phi)$
generating, by Stone's Theorem, $\Omega$ by
$$\Omega_{(t_1,\dots, t_n)}={\Delta_1}^{it_1}\dots {\Delta_n}^{it_n}.$$}\smallskip

\noindent{\bf Proof.} As $\omega$ is a $\phi$--invariant automorphism
group, it acts as a unitary group $\Omega$ on $L^2(\A, \phi)$.
$\B:=\B_{\phi,\omega}$ is the greatest subset of 
$I(\omega)\cap\N_\phi\cap{\N_\phi}^*$ invariant under the $\omega_z$,
$z\in{\Bbb C}$. $\B$ is a  $^*$--algebra invariant under the $\sigma_z$, $z\in{\Bbb C}$,
such that ${\cal B}^2$ is dense in $\B$ in the $2$--norm. 
We show that $\B$ is dense in $L^2(\A, \phi)$. For any element
  $x$ in $\N_\phi\cap{\N_\phi}^*$ (which is dense in $L^2(\A, \phi)$), 
$$x(\alpha):=(\alpha/\pi)^{n/2}\int_{{\Bbb R}^n}e^{-\alpha|t|^2}\omega_t(x)dt,$$ 
is still an element in $\N_\phi\cap{\N_\phi}^*$ which approximates $x$ as $\alpha\to+\infty$
 provided we show that $\Omega$ is strongly (or weakly) continuous.
One can show that, if $x\in\N_\phi\cap{\N_\phi}^*$
 then $\sigma_z(x(\alpha))\in\N_\phi\cap{\N_\phi}^*$, $z\in{\Bbb C}$, so that actually
$x(\alpha)\in\B$. We now show that $\Omega$ is weakly continuous. 
We first consider the case of the modular group $\sigma$.
 Since $\phi$ is a KMS weight,
 the function $t\to\phi(x^*\sigma_t(x))$ is continuous for $x\in\N_\phi$. In the 
general case, it is enough to show that the functions $\phi(\omega_t(x^*)y)$ 
are continuous when $x$, $y$ range over a subset of $\N_\phi$ dense in the 
$2$--norm. 
 Now $\omega$ is pointwise
norm continuous, therefore $\phi(x^*\omega_t(a)y)$ is continuous for $a\in\A$,
$x, y\in\N_\phi$. If in particular $a, x, y\in\cal{B}$, 
$$\phi(x^*\omega_t(a)y)=\phi(\sigma_{i/2}(\omega_t(a)y)\sigma_{i/2}(x)^*)=
\phi(\omega_t(\sigma_{i/2}(a))\sigma_{i/2}(y)\sigma_{-i/2}(x^*)),$$  
by the KMS property.
We can then consider the positive selfadjoint operators $\Delta_1,\dots,\Delta_n$
that generate the $n$--parameter group $\Omega$, by Stone's Theorem,
by $\Omega_{(t_1,\dots,t_n)}={\Delta_1}^{it_1}\dots{\Delta_n}^{it_n}$. 
Arguments similar to those previously used  show that if 
$x, y\in {\cal B}^2$, then $z\to\phi(x^*\omega_z(y))$ is an entire  function coinciding
with $\phi(\omega_{-z}(x^*)y)$ by the uniqueness principle of analytic functions.
${\cal B}$ is a common core for the $\Delta_j$'s.
It follows in particular that $\Omega_{z_1,\dots,z_n}$, regarded as an operator 
${\Delta_1}^{iz_1}\dots{\Delta_n}^{iz_n}$ with domain $\cal{B}$
is a densely defined preclosed operator. \smallskip

We shall consider the two 
 subspaces ${\mathcal B}_{\phi,\omega}\subset L^2(\A, \phi)$ and
${\mathcal B}_{\phi_r,\omega}\subset
L^2(\A, \phi_r)$ as 
natural common cores for certain unbounded  bijections
naturally arising from the locally compact quantum group.
For example, let us  write $S$ and $S_r$ for the closed 
operators $S_\phi$ and $S_{\phi_r}$ defined by the adjoint, 
denoting their polar decompositions by $S=J\Gamma^{1/2}$ 
and $S_r=J_r\Gamma^{1/2}_r$. $S_r$ and $S$ determine 
each other in the sense that there is an antiunitary operator,   
$Y: L^2(\A, \phi_r)\to L^2(\A, \phi)$ such that $Ya=R(a)^*$, $a\in\cal{A}$, 
intertwining them: $SY=YS_r$. Hence by the uniqueness of polar decompositions, 
we have $JY=YJ_r$, $\Gamma^{1/2}Y=
Y\Gamma_r^{1/2}$ as well. The next result describes other unbounded 
bijections arising from the coinverse. We denote by $\Delta$ and $\Delta_r$
the generators of the scaling group $t\mapsto\nu^{\frac{t}{2}}\tau_t$ when regarded as a 
strongly continuous unitary group on $L^2(\A, \phi)$ and $L^2(\A, \phi_r)$
respectively.
Clearly $\Delta Y=Y\Delta_r$. \medskip

\noindent{\bf 5.5 Theorem} {\sl Let $\sigma$ denote the modular group of 
$\phi$, $\tau$ the scaling automorphism group of the locally compact quantum 
group $\A$, and let $\omega$ be the $2$--parameter
group $\omega_{s,t}=\nu^{\frac{s}{2}}\tau_s\sigma_t$.
The following operators 
$${\mathcal B}_{\phi_r,\omega}\subset L^2(\A, \phi_r)\to {\mathcal 
B}_{\phi, \omega}\subset L^2(\A, \phi)$$
taking $a$ to $\kappa(a)$, $\kappa(a)^*$ and $\kappa(a^*)$ are densely 
defined and preclosed. Denoting the closures of the first 
operator by $K$, the closures of the second and third are $SK$ and $KS_r$, 
respectively. We have the following polar decompositions:
$$\nu^{+{\frac{i}{4}}}SK=Y\Delta_r^{1/2},$$
$$\nu^{-{\frac{i}{4}}}KS_r=Y\Delta_r^{-1/2},$$
$$\nu^{-{\frac{i}{4}}}K=JY\Gamma_r^{1/2}\Delta_r^{1/2}.$$ 
 Furthermore, 
$$(KS_r)^*=\nu^{i/2}(SK)^{-1},$$
$$(SK)^*=\nu^{-i/2}(KS_r)^{-1}.$$}\smallskip 

\noindent{\bf Proof.} 
It follows from the previous discussion that $\Delta_r^{-1/2}$, 
$\Delta_r^{1/2}$ and $\Gamma_r^{1/2}$ are bijective, densely defined, 
positive and  essentially selfadjoint on the indicated domains. 
In fact they are connected with $1$--parameter subgroups of $\omega_{s,t}$.
Since $Y$ is antiunitary, it is clear that $SK$ and $KS_r$
are bijections between the indicated domains.
 The polar decomposition of $SK$ and $KS_r$ on $\B_{\phi, \omega}$, 
 now follows
from the data of the locally compact quantum group. Furthermore $SK$ and $KS_r$ form part of an essentially
selfadjoint pair in the sense explained in the appendix.
 As $S$ and $S_r$ 
also form part of an essentially self--adjoint pair, it follows from Lemma
A.2 
that the operators in question can be denoted by $K$, $SK$ and $KS_r$. 
On the indicated domain one has:
 $$K=S(SK)=\nu^{{\frac{i}{4}}}J\Gamma^{1/2}Y\Delta_r^{1/2}=
\nu^{{\frac{i}{4}}}JY\Gamma_r^{1/2}\Delta_r^{1/2}.$$
Since $\B_{\phi, \omega}$ is $^*$--invariant and $\sigma_z$--invariant, $z\in
{\Bbb C}$,
one deduces, looking at the polar decomposition of $S$, that $J$ acts bijectively
on $\B_{\phi, \omega}$ too, and therefore 
 $K$ is a bijection from $\B_{\phi, \omega}$ to ${\mathcal B}_{\phi_r,
\omega}$ as well. 
It remains to show that $(KS_r)^*=\nu^{i/2}(SK)^{-1}$ and  $(SK)^*=\nu^{-i/2}(KS_r)^{-1}.$ 
We prove the latter relation, as the former follows by taking inverses and adjoints.
The polar decompositions show that $I(\omega)\cap\N_\phi$ is a 
core for both $(KS_r)^*$ and $\nu^{i/2}(SK)^{-1}$. But the two operators coincide on the 
core, completing the proof.\medskip

\noindent{\bf 5.6 Corollary}
{\sl $K$ is a closed intertwiner from $V_r$ to $V$, 
$$VI\otimes K= I\otimes K V_r.$$
In particular, if 
$U:=\nu^{{\frac{i}{4}}}JY : L^2(\A, \phi_r)\to L^2(\A, \phi)$  denotes
the polar part of $K$,
$$I\otimes U V_r=VI\otimes U.$$ 
Thus $V_r$ is a corepresentation of $V$. Moreover 
$$V(SK\otimes S)= (SK\otimes S)W,$$ 
where $W$ is the unitary on $L^2({\cal A},\phi_r)\otimes 
L^2({\cal A},\phi)$ 
defined by $W a\otimes b=\delta(b)a\otimes I$. In particular, taking the polar 
decomposition of $(SK\otimes S)^*$,
$$V\Delta^{1/2}\otimes\Gamma^{-1/2}=
\Delta^{1/2}\otimes\Gamma^{-1/2}V,$$
$$VY\otimes J=Y\otimes JW.$$ }\medskip

\noindent{\bf Proof.} 
We identify $K^*$ on a suitable $^*$--invariant core $\B$ of jointly
analytic 
vectors for $\sigma$ and $\tau$ contained in $\N_\phi$ (see Lemma 5.4). 
The polar decomposition of $K$ described in the previous Theorem shows
that a subset $\B\subset\A\subset L^2(\A,\phi)$ satisfies these properties
if $R(\B)^*=Y^*\B$
is a core for $K$  in $L^2(\A, \phi_r)$ satisfying similar properties. 
 Consider the right invariant weight $\phi_r=\phi\circ R$, with modular group
$x\mapsto\sigma^r_t(x)=R\sigma_{-t}(R(x))$. 
For $a\in\B$, $b\in R(\B)^*$,
$$(a, K(b))=\phi(a^*\kappa(b))=\nu^{i/2}\phi_r(b\kappa(a)^*)=$$
$$
\nu^{i/2}\phi_r(\sigma^r_{-i}(\kappa(a))^*b)$$
thus, for $a\in \B$, $K^*(a)=\nu^{-i/2}\sigma^r_{-i}(\kappa(a))$.
We need to show the following  relations
$$VI\otimes K< I\otimes K V_r,\eqno(5.4)$$
$$V_rI\otimes K^*< I\otimes K^* V.\eqno(5.5)$$
We start from (5.4). We claim that
it suffices to prove that, for $b\in R(\B)^*$,
$d\in\B$, $a,c\in\N_\phi$,
$$(Va\otimes\kappa(b), c\otimes d)=(V_ra\otimes b, c\otimes K^*
d).\eqno(5.6)$$
Indeed, the algebraic tensor product $\N_\phi\odot \B$ is a core
for $I\otimes K^*$, so  equation (5.6) shows that $V_r\N_\phi\odot
R(\B)^*$ lies in the domain of $(I\otimes K^*)^*=I\otimes K$ and
for $x\in\N_\phi\odot R(\B)^*$, $I\otimes K V_rx=VI\otimes Kx.$
On the other hand $\N_\phi\odot
R(\B)^*$ is a core for $I\otimes K$, so the claim follows.
We compute the left hand side of (5.6).
$$\phi\otimes\phi(a^*\otimes \kappa(b)^*\delta(d)c\otimes I)=
\phi(a^*\iota\otimes\phi(I\otimes\kappa(b)^*\delta(d))c)=$$
$$\phi(a^*\kappa(\iota\otimes\phi(\delta(\kappa(b)^*)I\otimes d))c)=
\phi(a^*\iota\otimes\phi\circ\kappa^{-1}(I\otimes \kappa(d)\vartheta
\circ\delta(b^*))c)=$$
$$\nu^{-i/2}\phi(a^*\iota\otimes\phi_r(I\otimes \kappa(d)
\vartheta\circ\delta(b^*))c)=
\nu^{-i/2}\phi(a^*\iota\otimes\phi_r(\vartheta\circ\delta(b^*)
I\otimes \sigma^r_{-i}(\kappa(d)))c)=$$
$$\nu^{-i/2}\phi\otimes\phi_r(a^*\otimes I\vartheta\circ\delta(b^*)
c\otimes \sigma^r_{-i}(\kappa(d))=
(V_ra\otimes b, c\otimes K^*d).$$
Arguing in a similar way, we see that $(5.5)$ will be a consequence of 
$$(Va\otimes b, c\otimes\kappa(d))=(V_ra\otimes K^*(b), c\otimes
d),\eqno(5.7)$$
for $a,c\in\N_\phi$, $b\in\B$, $d\in R(\B)^*.$
The r.h.s. equals
$$\nu^{i/2}\phi\otimes\phi_r(a^*\otimes
I\vartheta\circ\delta\circ\sigma^r_{i}\circ\kappa^{-1}(b^*))c\otimes d)$$
while the l.h.s. equals 
$$\phi\otimes\phi(a^*\otimes b^*\delta(\kappa(d))c\otimes I),$$
therefore it suffices to show that
$$\iota\otimes\phi(I\otimes
b^*\delta(\kappa(d)))=\nu^{i/2}\iota\otimes
\phi_r(\vartheta\circ\delta\circ\sigma^r_{i}\circ\kappa^{-1}(b^*)I\otimes
d).$$
Using
successively
$$\sigma^r_{i}\circ\kappa^{-1}=\kappa^{-1}\circ\sigma_{-i},$$
$$\vartheta\circ\delta\circ\kappa^{-1}=
\kappa^{-1}\otimes\kappa^{-1}\circ\delta$$
and
$$\nu^{i/2}\phi_r\kappa^{-1}=\phi$$
we see that the r.h.s. of the previous relation equals
$$\kappa^{-1}\iota\otimes\phi(I\otimes\kappa(d)\delta(\sigma_{-i}(b^*)))$$
which in turn equals
$$\iota\otimes\phi(\delta(\kappa(d))I\otimes\sigma_{-i}(b^*))$$
and the proof of (5.7) is completed by the KMS property of $\phi$.

We briefly sketch the second part of the proof.
We need to show the two relations
$$VSK\otimes S< SK\otimes S W,$$
$$W(SK)^*\otimes S^*<(SK)^*\otimes S^*V$$
which will follow respectively from 
$$(a\otimes b, VSK\otimes S c\otimes d)=(Wc\otimes d, (SK)^*\otimes
S^*a\otimes b),$$
for $a,b,d\in\B$, $c\in R(\B)^*,$
$$(a\otimes b, W(SK)^*\otimes S^* c\otimes d)=(Vc\otimes d, SK\otimes
S a\otimes b)$$
for $b,c,d\in\B$, $a\in R(\B)^*$. These equations
  can be obtained, in turn,
by computations similar to those of the first part of the proof.
\medskip

From the general structure of the quantum
groups under consideration, it follows that the dual Hopf algebra is equipped
with the same structure as $\A$ \cite{KvD}. 
We have already noted that the coinverse
$\kappa'$ of $\A_{\max}(V)$  serves to define the adjoint on $\A$
and consequently that   $\kappa'$  is uniquely determined by $\kappa$ via 
 $\kappa'(\omega)(a)=\omega(\kappa(a))$, $\omega\in{\hat\M}(V)_*$,
$a\in\A$, where $\A$ is being considered as a dense subspace of
$\hat\M(V)$.
Taking the square of $\kappa'$, which coincides with the square of the analytic
generator of $\tau'$, it follows that the coinverse data $R'$ and $\tau'$ of $\A_{\max}(V)$
are determined by those of $\A$ and similar formulas hold.  In particular,
$\kappa'(\omega)^*=\ov{\omega}$. We  may thus  write 
$$\ov{\omega}=\kappa'(\omega)^*=
{\tau'}_{i/2}(R'(\omega)^*), \quad\omega\in{\hat\M}(V)_*.$$
In fact, ${\tau'}_{i/2}$ is spatially implemented in the regular representation.
\smallskip

\noindent
{\bf 5.7 Proposition} {\sl Let $\tau'$ be 
the scaling group of $\A_{\max}(V)$ as defined above. Then 
for $\omega\in\D_{{\tau'}_{i/2}}$, $\xi'\in\D_{\Gamma^{-1/2}}$ and 
$\eta'\in\D_{\Gamma^{1/2}}$
$$(\xi',\pi_V({\tau'}_{i/2}(\omega))\eta')=(\Gamma^{-1/2}\xi',
\pi_V(\omega)\Gamma^{1/2}\eta').$$}\medskip

\noindent{\bf Proof.} 
It follows from Corollary 5.6 that 
$$\Delta^{it}\otimes\Gamma^{-it} V=V\Delta^{it}\otimes\Gamma^{-it},$$
Hence 
$$\tau_t\otimes i(V)=I\otimes \Gamma^{it}VI\otimes \Gamma^{-it},$$
Recalling that $\pi_V(\omega)=\omega\otimes i(V)$ and that $\tau'_t$ is 
the transpose of $\tau_t$, we deduce that
$$\pi_V(\tau'_t(\omega))=\Gamma^{it}\pi_V(\omega)\Gamma^{-it},\quad 
\omega\in\hat\M(V)_*.$$ 
Now with $\omega$, $\xi'$ and $\eta'$ as in the statement of the proposition, 
we have functions $z\mapsto(\xi',\pi_V(\tau'_z(\omega)\eta')$ and 
$z\mapsto(\Gamma^{-i\overline z}\xi',\pi_V(\omega)\Gamma^{-iz}\eta')$, 
analytic in the strip $0<\Im z<{1\over 2}$. Their boundary values agree 
on the real line and taking the boundary values at $z={i\over 2}$, we get
the required result.\medskip

  Now $j':=R'\circ *$ is a conjugation defined on $\A_{\max}(V)$ and 
satisfies $\delta\circ j'=j'\otimes j'\circ\theta\circ\delta$. Thus we 
know from Section 4 that $j'$ defines a conjugation on the tensor 
category of representations of $\A_{\max}$.\smallskip 

\noindent
{\bf 5.8 Corollary} {\sl A conjugate for $V$ in the conjugation 
defined by $j'$ is a conjugate for $V$ in the sense of Equation (5.1). 
The converse holds if $T$ in (5.1) has the form $T=J\Gamma^{1/2}$ 
with $J$ antiunitary.}\smallskip 

\noindent
{\bf Proof.} Let $J$ be an antiunitary operator then 
$\omega\mapsto J\pi_V(j'\omega)J^{-1}$, $\omega\in\hat\M(V)_*$, 
defines a representation of $\A_{\max}(V)$ and a conjugate $\overline V$ 
of $V$ in the conjugation defined by $j'$ with 
$$\omega\otimes i(\overline V)=J\pi_V(j'\omega)J^{-1}.$$ 
Now let $T:=J\Gamma^{1/2}$ and pick $\eta'\in\D_T$ and 
$\xi'\in\D_{T^*}$, i.e.\ $\eta',J\xi'\in\D_{\Gamma^{1/2}}$, then 
$$(\xi\otimes\xi',\overline V\eta\otimes T\eta')=
(\xi', \omega_{\xi,\eta}\otimes\iota(\overline V)T\eta')=$$ 
$$(\xi', J\pi_V(j'\omega_{\xi,\eta})J^{-1}T\eta')=(\xi',
J\pi_V(j'\omega_{\xi,\eta})\Gamma^{1/2}\eta')=$$
$$(\pi_V(j'\omega_{\xi,\eta})\Gamma^{1/2}\eta', \Gamma^{-1/2}T^*\xi').$$ 
 Now $j'\omega_{\xi,\eta}\in\D_{\tau_{i/2}}$, 
in fact
$\tau_{i/2}j'\omega_{\xi,\eta}=\overline\omega_{\xi,\eta}=\omega_{\eta,\xi}$.
Hence, by Proposition 5.7, 
$$(\pi_V(j'\omega_{\xi,\eta})\Gamma^{1/2}\eta',
\Gamma^{-1/2}T^*\xi')=(\pi_V(\tau'_{i/2}j'\omega_{\xi,\eta})\eta', T^*\xi')=$$
$$(\pi_V(\omega_{\eta,\xi})\eta', T^*\xi').$$
We conclude that 
$$(\xi\otimes\xi',\overline V\eta\otimes T\eta')=(V\xi\otimes\eta',\eta\otimes T^*\xi'),$$ 
so that $\overline V$ is a conjugate for $V$ in the sense of (5.1), as 
required. The converse follows by reversing the argument.\smallskip 

  Notice that, since $\Delta^{1/2}\Gamma^{-1/2}$ is a self--intertwiner of 
$V$, the result remains valid if we take $T$ to be of the form $T=J\Delta^{1/2}$ 
with $J$ antiunitary. We could replace $V$ by a corepresentation $W$ in the above 
argument if we knew that the automorphism group $\tau_t$ were unitarily 
implemented in the representation $\pi_W$.\smallskip

We now show, in analogy with the classical case, that the right regular 
corepresentation is a canonical choice of conjugate for the left regular one. 
We therefore look for a unitary corepresentation $\ov V$ and an antiunitary operator  $J$ such that 
$$(\xi\otimes \xi', \overline{V}\eta\otimes J\Delta^{1/2}\eta')=
(V\xi\otimes\eta', \eta\otimes\Delta^{1/2}J^*\xi'),$$
for $J^*\xi', \eta'\in\D_{\Delta^{1/2}}$, $\xi,\eta\in L^2(\A, \phi)$.\smallskip

\noindent{\bf 5.9 Theorem} {\sl The pair $\ov V= V_r$, 
$T=Y^*\Delta^{1/2}=\Delta_r^{1/2}Y^*$ solves the conjugate equation (5.1)
for 
the regular corepresentation $V$and shows that $V_r$ is a conjugate for $V$ 
in the conjugation defined by $j'$.}\medskip

\noindent{\bf Proof.} In view of the previous discussion, it suffices to 
verify the above equation and we begin by computing the l.h.s.\ $L$
making the change of variable $\hat\xi:=Y\Delta_r^{1/2}\xi'$ and using 
Theorem 5.5.
$$L=(\xi\otimes\Delta_r^{-1/2}Y^*\hat\xi,V_r\eta\otimes\Delta_r^{1/2}Y^*\eta')  
=(\xi\otimes(KS_r)^*\hat\xi, V_r\eta\otimes(KS_r)^{-1}\eta')=$$
$$=\nu^{-i/2}(\xi\otimes(SK)^{-1}\hat\xi, V_r\eta\otimes(KS_r)^{-1}\eta').$$
We choose $\hat\xi, \eta'$ belonging to a suitable common 
core of the indicated operators contained in $\N_\phi$ and write 
$(SK)^{-1}\hat\xi=\kappa(\hat\xi)^*$,
$(KS_r)^{-1}\eta'=\kappa({\eta'}^*)$. Then, recalling the 
definition of $V_r$, we have 
$$L=\nu^{-i/2}(\xi\otimes\kappa(\hat\xi)^*,V_r\eta\otimes\kappa({\eta'}^*))=
\nu^{-i/2}\phi\otimes\phi_r(\xi^*\otimes\kappa(\hat\xi)\theta\delta(\kappa(\eta^{'*}))
\eta\otimes I),$$
where the inner product refers to $L^2(\A, \phi)\otimes L^2(\A,\phi_r)$. Thus 
$$L=\nu^{-i/2}\phi(\xi^*\iota\otimes\phi_r(I\otimes\kappa(\hat\xi)
\kappa\otimes\kappa\circ\delta({\eta'}^*))\eta)=\nu^{-i/2}\phi(\xi^*
\iota\otimes\phi_r(\kappa\otimes\kappa(\delta(\eta^{'*})
I\otimes\hat\xi))\eta).$$ 
Since $\phi_r\circ\kappa=\nu^{i/2}\phi$ and 
$\kappa(\iota\otimes\phi(\delta(\eta^{'*})I\otimes\hat\xi))=
\iota\otimes\phi(I\otimes\eta^{'*}\delta(\hat\xi))$, we have  
$$L=\phi(\xi^*\kappa(\iota\otimes\phi(\delta({\eta'}^*)I\otimes\hat\xi))\eta)=
\phi(\xi^*\iota\otimes\phi(I\otimes\eta^{'*}\delta(\hat\xi))\eta)=$$
$$=(\xi\otimes\eta',V^*\eta\otimes\hat\xi)=(V\xi\otimes\eta',\eta\otimes\hat\xi),$$ 
completing the proof.
\medskip

We have seen above that we have 
 a conjugation $j':=R'\circ *$ defined on ${\cal A}_{\max}(V)$ 
and satisfying $\delta\circ j'=j'\otimes j'\circ\theta\circ\delta$,  
and we may, as described earlier, define a conjugation on its tensor
category 
of representations. We state this as a theorem.\medskip

\noindent{\bf 5.10 Theorem} {\sl Pick for each 
object $H$ in the image of 
the trivializing endofunctor $\iota$ on ${\cal C}(V)$ an antiunitary operator
$J_H: H\to \ov H$. Let $\pi_W$ be the representation of 
${\cal A}_{\max}(V)$ 
on $H$ corresponding to a $W\in{\cal C}(V)$. Let $\ov W$ denote the 
corepresentation of $V$ defined by the representation 
$A\mapsto J_H\pi_W(j'(A))J_H^{-1}$ of ${\cal A}_{\max}(V)$. 
Set 
$$J_W:=J_{\iota(W)},$$ 
whence
$$\ov{\iota(W)}=\iota(\ov{W}).$$
Given $T\in(W, W')$, we set 
$\ov T:=J_WT{J_{W}}^{-1}\in(\ov W,\ov{W'})$.
Then $T\in (W,W')\mapsto \ov T\in (\ov W,\ov W')$ is an antilinear 
functorial equivalence. We get a semilinear tensor $W^*$--category 
with conjugation with natural unitary equivalences 
$$d_W:=J_{\ov W}\circ J_W,$$ 
$$c_{W,W'}:=J_{W\times W'}\circ (J_W\times J_{W'})^*,$$ 
where $J_W\times J_{W'}:=\theta(\ov H,\ov H')\circ J_W\otimes
J_{W'}$.}\medskip
 
   Now that we have equipped the category of unitary corepresentations 
of a multiplicative unitary with a conjugation, it follows from Proposition 4.3 
that  the conjugate of a left regular object is a 
right regular object of the same category. 
\medskip

  On the other hand, we also see from Theorem 5.10 that the tensor 
$W^*$--category of finite--dimensional representations of a compact quantum  
group admits two canonical conjugations, which are in general quite 
distinct. The one coming from   Theorem 5.10 is embedded in a semilinear 
tensor $W^*$--category of Hilbert spaces and is related to the antilinear 
involution $j'$ which commutes with the adjoint. The other cannot be 
embedded unless the intrinsic dimensions are integral (Corollary 4.2) 
and is related to the antilinear involution $^*\circ\kappa'$ which does 
not commute with the adjoint in general. In the case of $SU_q(2)$, it 
is clear that the conjugation on objects is the same in both cases 
since the fusion rules imply that each irreducible is self--conjugate. We 
shall show that this holds for all compact quantum groups. 
Now $^*\circ\kappa'=\tau'_{i/2}\circ j'=j'\circ\tau'_{-i/2}$, so 
the relation between the two conjugations should be describable in terms 
of $\tau'$. To this end, we consider the category of finite--dimensional 
unitary representations of the compact quantum group. As is well known, 
there is an associated Hopf--von Neumann algebra. Its elements are most 
conveniently described as bounded natural transformations of the 
embedding functor $F$ into the underlying Hilbert spaces. Now we have 
identified a natural transformation $\rho\mapsto f_\rho$ prior to 
Theorem 4.1. In general $f$ is not an element of the Hopf--von Neumann 
algebra $(F,F)$ since the natural transformation is not bounded in 
general. Nevertheless, its bounded functions such as $\rho\mapsto f_\rho^{it}$ 
will be elements of $(F,F)$ and these induce the inner automorphisms 
$\tau'_t$ of $(F,F)$. It follows that $f$ itself induces $\tau'_{-i}$. 
The relation between the two conjugations is now described in the 
following proposition.\smallskip 

\noindent
{\bf 5.11 Proposition } {\sl Adjoin to the tensor $W^*$--category of
finite 
dimensional representations of a compact quantum group two sets of antilinear 
intertwiners: $X\in(\rho,\sigma)_\kappa$ defined by  
$$X\rho(^*\circ\kappa'(A))=\sigma(A)X,\quad X\in\D_{^*\circ\kappa},$$ 
and $Y\in(\rho,\sigma)_j$ defined by 
$$Y\rho(j'(A))=\sigma(A)Y,\quad A\in(F,F).$$ 
Then we obtain two conjugations ${\cal T}_\kappa$ and ${\cal T}_j$. 
The second has the induced $^*$--structure, whilst the first has 
the $^*$--structure  ${\cal T}^{\ov a}$ described in
Theorem 4.1 b). There is an isomorphism 
of tensor categories which is the identity on linear arrows and 
takes $X\in(\rho,\sigma)_\kappa$ into $Xf_\rho^{-1/2}=f_\sigma^{1/2}X$.
This isomorphism does not commute with the adjoint in general.}\smallskip 

\noindent
{\bf Proof.} ${\cal T}_j$ is the conjugation defined by $j'$ and has the 
induced $^*$--structure. To see that ${\cal T}_\kappa$ is the conjugation
${\cal T}^{\ov a}$ 
described in Theorem 4.1 b), it suffices to show that the invertible 
antilinear operators ${\ov T}_\rho$, introduced in connection with
Theorem 4.1 
are in $(\rho,\ov\rho)_\kappa$. Now this has already been
noted after $(5.3)$.  A
computation shows that 
$X\in(\rho,\sigma)_\kappa$ if and only if $Xf_\rho^{-1/2}\in(\rho,
\sigma)_j$, so we 
will have the desired isomorphism once we show that 
$Xf_\rho^{-1/2}=f_\sigma^{1/2}X$. However $T_\rho^*T_\rho=f_\rho$ and 
$T_\rho^{-1*}T_\rho^{-1}=f_{\ov\rho}$. Hence 
$f_{\ov\rho}T_\rho=T_\rho^{-1*}=T_\rho f_\rho^{-1}$, so 
$f_{\ov\rho}\ov{T}_\rho=\ov{T}_\rho {f_\rho}^{-1}$ as required. From
this
we deduce that if $X\in(\pi,\rho)_\kappa$, then 
$\ov{T}_\rho\circ X\in(\pi,\ov\rho)$ so 
$$\ov{T}_\rho Xf_\pi^{-1}=f_{\ov\rho}^{-1}\ov{T}_\rho X=\ov{T}_\rho
f_\rho X,$$ 
leading to the desired result.\smallskip

We see, therefore, that for compact quantum groups the conjugations
$\T_\kappa$ and $\T_j$ although corresponding to different notions
of antilinear intertwiner lead to the same notions of conjugate
object.\smallskip

We conclude this section asking, 
in view of Theorem 2.3 and Proposition 4.3, whether there is a
relationship
between conjugation and standard braided symmetries for $\T$.
Assume that
$\T$ admits a conjugation $J$
assigning to a left regular object $V$  a conjugate  $\ov V$
equivalent to $V$ itself.
By Proposition 4.3 $\ov V$ is a right regular object, so $V$  is a right
regular object of $\T$ as well. Explicitly, if $U\in(\ov V, V)$
is a unitary and $\eta^{\ov V}\in(R^{\ov V}, R^{\ov V}\iota)$ is the natural 
unitary transformation derived from $J$ as in the proof of Proposition 
4.3,
then
$$\eta_W:=1_{\iota(W)}\times U\circ{\eta^{\ov V}}_W\circ 1_W\times U^*=$$
$$1_{\iota(W)}\times U\circ J_V\times 
{J_{\iota(W)}}^{-1}\circ\xi_{\ov W}\circ J_W\times{J_V}^{-1}\circ 1_W\times U^*$$
makes $V$ into a right regular object.
To have a standard braided symmetry we further need the commutation relation  
d') of Lemma 2.4
$$\eta_V\times1_{\iota(V)}\circ 1_V\times\xi_V=1_V\times\xi_V\circ
\eta_V\times1_{\iota(V)},$$
which
 expresses  a precise  relationship between
 the antilinear conjugation arrows $J_V$, $J_{\iota(V)}$, the unitary 
intertwiner
$U\in(\ov V, V)$        and $\xi_V$.
We next interpret this relation from the viewpoint of the Banach Hopf
algebra
associated to $V$.\medskip

\noindent{\bf 5.12 Proposition}
{\sl Let $\T$ be a tensor $C^*$--category
with a left regular object $V$ and conjugation  as described in
Proposition 4.3. Assume furthermore that
there is a unitary $U\in(\ov V, V)$  such that
$UJ_V\hat\A(V^*){(UJ_V)}^{-1}$ and $\hat{\A}(V)$ commute. 
Then $\T$ has a standard braided symmetry.}\smallskip

\noindent{\bf Proof} Let as usual $\A(X)$ and $\hat\A(X)$ stand for the 
Banach spaces derived from a unitary operator $X$ on a tensor product
 Hilbert space
compressing its first and the second component respectively.
First note that, $V$ and $\ov V$ being unitarily equivalent, 
$\hat{\A}(\ov V)=\hat\A(V)$,
so 
$$\A({\eta^{\ov V}}_V)=J_V\hat\A({\ov V}^*){J_V}^{-1}=J_V\hat\A(V^*){J_V}^{-1},$$
therefore 
$$\A(\eta_V)=U\A({\eta^{\ov V}}_V)U^*=UJ_V\hat\A(V^*)(UJ_V)^{-1}.$$
Recall now that  d$'$) of Lemma 2.4 is equivalent to requiring that
$\A(\eta_V)$ and $\hat\A(V)$ commute.\medskip

In the case where $J$ is a conjugation taking $V$ to a conjugate $\ov V$ solving the conjugate
equation,  we have already noted that $\hat\A(V)=\hat\A(V^*),$
therefore one is reduced to requiring that 
$$UJ_V\hat\A(V){(UJ_V)}^{-1}\subset\hat\A(V)'.$$
We now show that for locally compact quantum groups there are canonical
choices of $U$ and $J_V$ for which the above commutation relation holds.
Indeed in this case a solution
of the conjugate equation for the regular corepresentation
is given by $J_V=Y^*$ and  $\ov V=V_r$ (Theorem 5.9). Furthermore
the polar part $U=JY$ of $K$ is a unitary intertwiner
in $(\ov V, V)$ (Corollary 5.6), therefore $UJ_V=J$ is the polar
part of the Tomita operator on $L^2(\A, \phi)$, which
takes, via its adjoint action, $\hat\A(V)=\A$
into its commutant. We can thus conclude that the corepresentation 
category of $V$ has a standard braided symmetry.
Making use of the previous observation, we compute explicitly the
associated natural unitary transformation 
$\eta$ making $V$ into a right regular object.
First we have that 
$$\eta^{\ov V}_V=J_V\times {J_{\iota(V)}}^*\xi_{\ov V}J_V\times
{J_V}^*=$$
$$\vartheta_{H,K} Y^*\otimes Y {V_r}^*Y\otimes
Y^*\vartheta_{K,H}=\vartheta_{H, K}
W^*\vartheta_{K, H}$$
where $H=L^2(\A, \phi_r)$ and $W$ is the operator defined in Corollary
5.6. 
The last equality follows from the commutation relations between
$V$, $V_r$ and $W$ obtained in that corollary.
We next have that $\eta$  
is determined  by
$$\eta_V=1_{\iota(V)}\times U\circ\eta^{\ov V}_V\circ 1_V\times U^*=$$
$$I\otimes U\vartheta_{H, K} W^*\vartheta_{K, H} I\otimes U^*$$
where $U$ is the polar part of $K$. 
This equation may be understood as a categorical interpretation
in terms of the conjugation structure in $\C(V)$ of the usual way of
getting 
standard braided symmetries for quantum groups described in 
Proposition 2.11. (Unitaries of the form of $\eta_V$ had previously
appeared in
\cite{BS}
in
the
context of irreducible multiplicative unitaries). We have thus shown the
following result.\medskip

\noindent{\bf 5.13 Proposition} {\sl Let $(\A, \delta, R, \tau, \phi)$
be a locally compact quantum group, and let us endow
$\C(V)$ with a conjugation as described in Theorem 5.10
and let $V$ be the usual multiplicative unitary
associated to it. Then there is an
associated 
standard
braided 
symmetry $\varepsilon$ on $\C(V)$ whose evaluation in $V$
is given by
$$\varepsilon_V=V U\times I W^* U^*\times I\vartheta,$$}\medskip
where $W$ is the unitary on $L^2(\A, \phi_r)\otimes L^2(\A, \phi)$
defined in Corollary 5.6.

\section{Appendix. Essentially self--adjoint pairs}

The basic difficulty in defining conjugation in the 
context of quantum groups or multiplicative unitaries is that one starts from 
an antilinear involution which is not, in general, a conjugation on 
a $C^*$--algebra in that it may not commute with the adjoint and 
may only be densely defined and unbounded. In the theory of 
Woronowicz \cite{W}, \cite{WTK} 
this involution defines the contragredient representation and he shows 
how to pass to the conjugate representation in the context of compact 
quantum groups. This involves a variant of modular theory and, to prepare for 
this, we begin with a simple result on densely defined semilinear 
bijections, where semilinear is understood to mean linear or antilinear.
If $s$ is a densely defined semilinear mapping between Hilbert spaces 
$\cal{H}$ and $\cal{K}$ and $f$ a densely defined semilinear 
mapping from $\cal{K}$ to $\cal{H}$ then we call the pair $s$, $f$ 
Hermitian, essentially self--adjoint or selfadjoint if the matrix 
$$\left(\begin{matrix} 0&f \cr s&0\end{matrix}\right)$$ 
defines a semilinear mapping on ${\cal H}\oplus{\cal K}$ with the 
corresponding property.\smallskip 

\noindent
{\bf A.1 Lemma} {\sl Let $s$, $f$ be an essentially selfadjoint pair of 
semilinear bijections between dense subspaces of Hilbert spaces $\cal{H}$ 
and $\cal{K}$. Let $\hat s$ and $\hat f$ denote the semilinear involutions 
given by the matrices 
$$\left(\begin{matrix} 0&s^{-1}\cr s&0\end{matrix}\right),\quad 
\left(\begin{matrix} 0&f\cr f^{-1}&0\end{matrix}\right).$$ 
Then $\hat s\subset\hat f^*$.
The eigenspaces of $\hat s$ and $\hat f$ corresponding to the eigenvalue $\pm 1$ 
consist of vectors of the form 
$$\left(\begin{matrix} h\cr \pm sh\end{matrix}\right),\,\, \text{and}\,\, 
\left(\begin{matrix} h'\cr \pm f^{-1}h'\end{matrix}\right),$$ 
for $h\in{\cal D}_s$ and $h'\in{\cal D}_{f^{-1}}$.
Let the  closure of these eigenspaces be denoted $M^s_\pm$ and $M^f_\pm$ then 
$M^f_\pm=M^{s\perp}_\mp$. Writing $\hat S$ for the closure of $\hat s$ and 
denoting its polar decomposition by
$$\hat S=\hat J\hat\Delta^{1/2},$$ 
$\hat S$ and $\hat J$ are involutions, $\hat J\Delta^{1/2}=\Delta^{-1/2}\hat J$ 
and $\hat J(M^s_\pm)=M^f_\mp$.  
Furthermore $\hat S^*=\hat F$, the closure of $\hat f$.}\smallskip

\noindent
{\bf Proof} $\hat s$ and $\hat f$ are obviously involutions 
with eigenspaces as stated and $\hat f\subset\hat s^*$. The closures 
of these eigenspaces are obviously the eigenspaces of the closures 
$\hat S$ and $\hat F$. The eigenspaces of the involution $\hat S^*$ are 
then $M_-^{s\perp}$ and $M_+^{s\perp}$. But $\hat S^*=\hat F$ since the pair 
$s$, $f$ is essentially selfadjoint so $M^f_\pm=M^{s\perp}_\mp$. Finally, 
$\hat J(M^s_\pm)=M^f_\mp$ since $\hat S\hat J=\hat J\hat S$ and this and 
the remaining assertions just follow from the fact that $\hat S$ is a closed 
densely defined involution. 
\smallskip 

  The formalism here is best known in the antilinear case, through 
modular theory. In this case, it is usual to formulate the result on the 
eigensubspaces by saying that $M^s_+$ is a standard subspace and $M^f_+$ 
is its symplectic complement. The above formulation brings out the similarities 
between the linear and antilinear case. We draw the obvious conclusions 
about the original operators thus 
$$\hat S=\left(\begin{matrix} 0&S^{-1}\cr S&0\end{matrix}\right),  \quad 
\hat\Delta=\left(\begin{matrix} S^*S&0\cr 0&S^{-1*}S^{-1}\end{matrix}\right).$$
Finally, letting $J\Delta^{1/2}$ be the polar decomposition of $S$, 
$$\hat J=\left(\begin{matrix} 0&J^{-1}\cr J&0\end{matrix}\right).$$

  To generalize the above version of modular theory so as to treat 
compositions of mappings,  we let ${\cal V}$ be a subcategory of the category of 
vector spaces with linear (or antilinear) maps, and pick, for each object
$V$ of ${\cal V}$, an object $V'$ of ${\cal V}$ such that 
$V$ and $V'$ are vector spaces in 
duality, via a bilinear or sesquilinear form. Let
$s\in(V_1, V_2)$ be an   arrow endowed with a transposed
arrow $f=s'\in({V_2}', {V_1}')$ with respect to the duality. For example, if
$s$ is linear,
$f$ is the linear map defined by
$$({k_2}', sk_1)=(f{k_2}', k_1),\quad k_1\in V_1, {k_2}'\in{V_2}'.$$
While, if $s$ is antilinear, $f$ is the antilinear  map defined by
$$({k_2}', sk_1)=\overline{(f{k_2}', k_1)},\quad k_1\in V_1, {k_2}'\in{V_2}'.$$

This would seem to be a natural notion to study in the theory of dual 
spaces. In the present context we need a more restrictive notion, replacing 
vector spaces by (positive--definite) scalar product spaces and requiring the 
duality to make the scalar product space $V'$ a dense subspace of the 
completion of $V$.
A pair of the form $(s, f)$ will in this case  be called 
an Hermitian pair of ${\cal V}$.
A particular 
instance of this situation would be to fix a set of Hilbert spaces, 
each Hilbert space $H$ having two distinguished dense subspaces
 $H^s$ and $H^f$, 
say, and take an arrow in the category to be a linear (or just semilinear) 
mapping from some $H^s_1$ into some $H^s_2$ whose adjoint is defined on 
$H^f_2$ and maps it into $H^f_1$.\smallskip 

\noindent
{\bf A.2 Lemma} {\sl Let  ${\cal SP}^a$ be the category
of scalar product spaces 
with semilinear mappings. Then 
 for each hermitian pair $(s,f)$ of ${\cal SP}^a$, we  let
$\ov{s}$ and $\ov{f}$ denote the closure of $s$ and $f$, respectively.
If  $s's$ defines an  essentially self--adjoint pair, then
$\overline{s's}$ is the closure of $\ov{s'}\ov{s}$,
 and $\overline{ff'}$ the closure of $\ov{f}$ $\ov{f'}$.}\smallskip 

\noindent
{\bf Proof.} Let $h\in {\cal D}_{\ov{s}}$, $\ov{s}h\in{\cal
D}_{\ov{s'}}$ and 
$k$ in the domain of $f'$, then 
$$(k,\ov{s'}\ov{s}h)=(f'k,\ov{s}h)=(ff'k,h),$$ 
since $f'\subset {s'}^*$, $f\subset s^*$. Hence 
$s's\subset \ov{s'}\ov{s}\subset (ff')^*=\overline{s's}$ since $s's$ and 
$ff'$ 
is an essentially self--adjoint pair, completing the 
proof.\smallskip

  An instructive example of the situation of Lemma A.1 is got by 
considering a von Neumann algebra $\cal{M}$ with cyclic and separating vector 
$\Omega$. Then $\cal{M}$ and $\cal{M}'$ are in (sesquilinear) duality via 
$(M'\Omega,M\Omega)$. We may now take $\cal{C}$ to be the category of 
semilinear mappings between $\cal{M}$ and $\cal{M}'$ with dual. For 
example if $\kappa$ is a linear mapping from $\cal{M}$ to $\cal{M}$ then 
$\kappa$ has a dual 
if there is a linear mapping $\kappa'$ from $\cal{M}'$ to $\cal{M}'$ such that 
$$(\kappa'(M')\Omega,M\Omega)=(M'\Omega,\kappa(M)\Omega),\quad M\in\cal{M}\,\,,
M'\in\cal{M}',$$ 
and a linear mapping $\lambda$ from $\cal{M}$ to $\cal{M}'$ has a dual if 
there is a linear mapping $\lambda'$ from $\cal{M}$ to $\cal{M}'$ 
such that 
$$(\lambda'(M)\Omega,N\Omega)=(M\Omega,\lambda(N)\Omega),\quad M,\,N\in\cal{M}.$$
The adjoint on $\cal{M}$ has, of course, as its dual the adjoint on $\cal{M}'$. 
The modular conjugation $j$ from $\cal{M}$ to $\cal{M}'$ is selfdual. The 
adjoint of the modular automorphism $\sigma_t$ of $\cal{M}$ is the modular 
automorphism $\sigma'_t$ of $\cal{M}'$. In each case, we get an 
essentially self--adjoint pair.\smallskip

  In the above example, our scalar product spaces are actually left Hilbert 
algebras and in this case it is interesting to note what happens when 
the mappings are multiplicative and bijective (without necessarily 
commuting with the adjoint). The basic information is presented in the form 
of a lemma without proof.\smallskip

\noindent
{\bf A.3 Lemma } {\sl Let $s$, $f$ be an Hermitian pair of semilinear 
multiplicative bijections of left Hilbert algebras $\cal{A}$ and $\cal{A}'$ 
then for the associated closed operators $S$ and $F$ we have 
$$Sa=s(a)S,\quad S^*a=s^*(a)S^*,\quad a\in\cal{A},$$ 
$$Fa'=f(a')F,\quad F^*a'=f^*(a')F^*,\quad a'\in\cal{A}',$$ 
where $s^*(a)=s^{-1}(a^*)^*$, $f^*(a')=f^{-1}(a^*)^*$.}\smallskip

  Among the many variations of the above example, we can, in the above 
choose a faithful semifinite normal weight $\phi$ on $\cal{M}$ with associated 
left ideal ${\cal N}_\phi$ and let ${\cal B}_\phi$ denote the weakly dense 
$^*$--subalgebra of $\cal{M}$ consisting of entire elements $B$ for the modular 
automorphism group $\sigma^\phi$ such that 
$\sigma^\phi_z(B)\in{\cal N}_\phi
\cap{\cal N}_\phi^*$ for all $z\in{\mathbb C}$. Let ${\cal B}'_{\phi'}$ be the 
corresponding subalgebra of $\cal{M}'$ defined using the weight 
$\phi':=\phi\circ j^{-1}$ and consider the duality between ${\cal B}_\phi$ 
and ${\cal B}'_{\phi'}$.

  We are led to a simple example of the situation in Lemma A.2. 
Consider the situation of Lemma 5.4. We let 
${\cal C}_{\phi,\omega}$ be the set of semilinear mappings of 
${\cal B}_{\phi,\omega}$ whose adjoint is defined on 
${\cal B}_{\phi,\omega}$ 
and maps it into itself.
Elements of $\cal{C}_{\phi,\omega}$ include the restrictions of the
$\omega_z$, 
$z\in{\mathbb C}$ as well as the restriction of $S_\phi$. Taking the
closures, we recover the $\omega_z$ and $S_\phi$.\smallskip

 We give a rather more complicated example involving two
objects ${\mathcal B}_{\phi, \omega}$ and
${\mathcal B}_{\ov{\phi},\omega}$ motivated by Theorem 5.5.\smallskip 

\noindent
{\bf A.4 Lemma} {\sl Let $\phi$ be a lower semicontinuous densely defined 
weight on a $C^*$--algebra $\A$ equipped  with a conjugation $j$. Let 
$\omega: {\mathbb R}^n\to\Aut(\A)$ be a pointwise norm continuous 
$\phi$--invariant automorphism group of $\A$ and suppose that 
$j\omega_{{\mathbb R}^n}j=\omega_{{\mathbb R}^n}$. Let $\ov\phi:=\phi\circ
j$ 
then $j({\cal B}_{\phi,\omega})={\cal B}_{\ov\phi,\omega}$. Now let 
$\cal{C}$ denote the category of semilinear mappings between the two 
objects ${\cal B}_{\phi,\omega}$ and ${\cal B}_{\ov\phi,\omega}$ whose 
adjoints restrict to a mapping in the other direction. Then 
$\cal{C}$ contains ${\cal C}_{\phi,\omega}$ and ${\cal
C}_{\ov\phi,\omega}$ 
as categories in a natural way and the map $B\mapsto j(B)$, 
$B\in{\cal B}_{\phi,\omega}$ is an antilinear bijection.}\smallskip
 
\noindent
{\bf Proof.} We have $j({\cal N}_\phi)={\cal N}_{\ov\phi}$ and since $j$ 
normalizes $\omega_{{\mathbb  R}^n}$, $j({\cal B}_{\phi,\omega})=
{\cal B}_{\ov\phi,\omega}$.
 The remaining statements are now obvious.\smallskip 

Theorem 5.5 corresponds to taking for $\phi$ the left Haar weight
on a locally compact quantum group and a conjugation $j=R\circ^*$.

\end{document}